\documentclass[12pt, letterpaper]{article}
\usepackage[left=3cm, right=3cm, top=2cm, bottom=2cm]{geometry}
\pdfoutput=1
\usepackage{hyperref}
\hypersetup{colorlinks,allcolors=blue}
\usepackage[utf8]{inputenc}
\usepackage[backend=biber, sorting=none, autolang=hyphen,
natbib=true,maxbibnames=5,maxcitenames=2,uniquelist=false,url=false,doi=false,isbn=false,eprint=false]{biblatex}
\usepackage[english]{babel}
\usepackage{csquotes}
\usepackage{amsfonts}
\usepackage{overpic}
\usepackage[justification=centering,hypcap=true]{caption}
\usepackage{hypcap}
\usepackage{iftex,amsmath}
\usepackage{tgtermes}
\usepackage{float}
\usepackage{amssymb}
\usepackage{bbm}
\usepackage{faktor}
\usepackage[table]{xcolor}
\usepackage{array}
\usepackage[nomarkers,nolists]{endfloat}
\usepackage{tikz}
\usetikzlibrary{angles,patterns,calc}
\usetikzlibrary{arrows}
\usepackage{crossreftools}
\usepackage{wrapfig}
\usepackage{subfig}
\usepackage{makecell}
\counterwithin*{figure}{section}
\counterwithin*{figure}{subsection}
\counterwithin*{figure}{subsubsection}
\counterwithin*{figure}{paragraph}
\newcommand{\notifzero}[1]{%
  \ifnum\value{#1}>0
    \arabic{#1}.%
  \fi
}

\newcommand{\bb}[1]{\mathbb{#1}}
\newcounter{figures}[subsection]

\newcounter{figuress}[subsection]

\newcounter{allcounter}[subsection]
\newenvironment{prop}[1][]{\refstepcounter{allcounter}\par\medskip \textbf{Proposition~\thesubsection.\theallcounter \,\, \rmfamily}}{\medskip}

\newenvironment{deff}[1][]{\refstepcounter{allcounter}\par\medskip \textbf{Definition~\thesubsection.\theallcounter \,\, \rmfamily}}{\medskip}

\newenvironment{notation}[1][]{\refstepcounter{allcounter}\par\medskip \textbf{Notation~\thesubsection.\theallcounter \,\, \rmfamily}}{\medskip}

\newenvironment{theorem}[1][]{\refstepcounter{allcounter}\par\medskip \textbf{Theorem~\thesubsection.\theallcounter \,\, \rmfamily}}{\medskip}

\newenvironment{remark}[1][]{\refstepcounter{allcounter}\par\medskip \textbf{Remark~\thesubsection.\theallcounter \,\, \rmfamily}}{\medskip}

\newenvironment{conc}[1][]{\refstepcounter{allcounter}\par\medskip \textbf{Conclusion~\thesubsection.\theallcounter \,\, \rmfamily}}{\medskip}

\newenvironment{lemma}[1][]{\refstepcounter{allcounter}\par\medskip \textbf{Lemma~\thesubsection.\theallcounter \,\, \rmfamily}}{\medskip}

\newenvironment{note}[1][]{\refstepcounter{allcounter}\par\medskip \textbf{Note~\thesubsection.\theallcounter \,\, \rmfamily}}{\medskip}

\newenvironment{theoremm}[1][]{\refstepcounter{allcounter}\par\medskip \textbf{Theorem~\thesection.\theallcounter \,\, \rmfamily}}{\medskip}

\newenvironment{remarkk}[1][]{\refstepcounter{allcounter}\par\medskip \textbf{Remark~\thesection.\theallcounter \,\, \rmfamily}}{\medskip}

\makeatletter
\newcommand*{\rom}[1]{\expandafter\@slowromancap\romannumeral #1@}
\makeatother

\newcommand{\proof}{\emph{Proof.}\,\,\,\,}

\DeclareDocumentCommand{\inp}{g}
{
	\IfNoValueF{#1}{\langle \, #1 \, \rangle}
	\IfNoValueT{#1}{\langle \, \cdot\, ,\, \cdot \, \rangle}}

\newcommand{\pare}[1]{\left({#1}\right)}

\newcommand{\brac}[1]{\left\{ {#1}\right\}}

\newcommand{\mat}[1]{\pare{\begin{matrix} #1 \end{matrix}}}

\DeclareDocumentCommand{\dv}{d[]om}{
	\IfNoValueF{#1}{
		\IfNoValueF{#2}{
			\frac{{\rm d}^{#1}#2}{{\rm d}{#3}^{#1}}
		}
		\IfNoValueT{#2}{
			\frac{{\rm d}^{#1}}{{\rm d}{#3}^{#1}}
		}
	}
	\IfNoValueT{#1}{
		\IfNoValueF{#2}{
			\frac{{\rm d}#2}{{\rm d}#3}
		}
		\IfNoValueT{#2}{
			\frac{{\rm d}}{{\rm d}#3}
		}
	}
}
\DeclareDocumentCommand{\dd}{o}{
	\IfNoValueF{#1}{\,{\rm d}#1}
	\IfNoValueT{#1}{\,{\rm d}}
}

\DeclareDocumentCommand{\H}{}{\mathbb{H}}

\DeclareDocumentCommand{\N}{}{\mathbb{N}}
\DeclareDocumentCommand{\Q}{}{\mathbb{Q}}
\DeclareDocumentCommand{\R}{}{\mathbb{R}}
\DeclareDocumentCommand{\Z}{}{\mathbb{Z}}

\DeclareDocumentCommand{\cF}{}{\mathcal{F}}

\DeclareDocumentCommand{\cH}{}{\mathcal{H}}

\newcounter{lecture}[subsection]

\newcounter{classexercise}[subsection]

\newcounter{ex}[subsection]

\newcounter{allcontent}[subsection]

\newcounter{examples}[subsection]

\newcounter{claim}[subsection]

\newcounter{corollary}[subsection]

\DeclareDocumentCommand{\dImage}{r[] m m d[]}{
\begin{figure}[h]
\label{#3}
\begin{center}
\includegraphics[width=#1\textwidth]{#2}
	\IfNoValueF{#4}{\caption{#4}}
\end{center}
\end{figure}
}

\definecolor{lightapricot}{rgb}{0.99, 0.84, 0.69}
\definecolor{mossgreen}{rgb}{0.62, 0.87, 0.68}
\definecolor{grannysmithapple}{rgb}{0.66, 0.89, 0.63}
\definecolor{teagreen}{rgb}{0.82, 0.94, 0.75}
\definecolor{lightpastelpurple}{rgb}{0.69, 0.61, 0.85}
\definecolor{lilac}{rgb}{0.78, 0.64, 0.78}
\definecolor{pastelyellow}{rgb}{0.99, 0.99, 0.59}
\definecolor{palepink}{rgb}{0.98, 0.85, 0.87}
\definecolor{lightcyan}{rgb}{0.88, 1.0, 1.0}

\usepackage{eucal}

\title{Studying knots in self-covers of the modular flow}
\author{Sivan Eldar and Stav Fahima}
\date{}

\begin{document}

\maketitle

\textbf{Abstract.} Geodesic flows are a fascinating subject of study in dynamical systems, which was contributed to by many mathematicians. In this paper, we provide a combinatorial tool to help study the topological properties of modular knots. We construct templates for the infinitely many Anosov flows on the trefoil complement, which are finite order lifts of the geodesic flow on the modular surface. This construction is obtained by lifting Ghys’ modular template \cite{Ghys:07} using the covers as in \cite{clay:20}. This allows us to study the knot properties of periodic orbits in these flows, and to explicitly construct an infinite family of links of two components, one being the trefoil, all commensurable to one another.
\bigskip

\section{Introduction}

The set of two-dimensional lattices up to scaling is a topological space originally studied by Gauss in the context of number theory. There is an embedding in $S^3$ of the space of lattices that is homeomorphic with the trefoil complement \cite{mil71}. Considering the action of the matrices $\brac{\text{diag}(e^t,e^{-t})}_{t\in \R}$ on the space of lattices, every matrix in $SL_2(\Z)$ acting on a lattice returns the same lattice up to a change of base. From this action, we obtain a dynamical system which is the geodesic flow on the modular surface, a well-known example of an Anosov flow. This space has been studied thoroughly in the context of dynamical systems \cite{dist12},\cite{wen08}.
\medskip

Birman and Williams \cite{birman1:83},\cite{birman2:83} proved that every Anosov flow has a combinatorial representation as a template, that is, a branched surface that carries the periodic orbits of the flow up to isotopy and enables the study of their knot properties. They focused on the Lorenz template and showed that it determines an infinite family of knots and links, which they called \emph{Lorenz links}.
\medskip

Ghys \cite{Ghys:07} constructed a template, shown in figure 1.1 (a), for the geodesic flow on the modular surface which he named \emph{the modular flow}. He was the first to study the topological knot properties of its periodic orbits, called \emph{modular knots}, by proving that modular knots coincide with Lorenz knots. That is, every modular knot is isotopic in $S^3$ to one of the periodic orbits of the Lorenz attractor. He used the modular template to prove the equivalency between the Rademacher function of a periodic orbit and the orbit's linking number with the trefoil, giving it a topological meaning. He also proved that these knots are prime and fibered when embedded in $S^3$, relying on the work of Birman and Williams. His work inspired further topological results on these knots.
\medskip

The trefoil complement in $S^3$, on which the modular flow is defined, fibers over $S^1$. Thus, it admits an $n$-sheeted self-cover $p_n: S^3\setminus T \longrightarrow S^3\setminus T$, for $n=6k+1$ with $k \in \N_{>0}$, and $T$ being the trefoil (see \cite{rolfsen:03} for example). 
Given a geodesic flow on the trefoil complement, a significant result by Clay and Pinsky \cite{clay:20} shows that the pre-images of this flow under the covering maps $p_n$ yield non-equivalent Anosov flows on the trefoil complement for different values of $n$, hence proving it is possible to construct arbitrarily many such flows. In this paper, we extend their work by providing explicit combinatorial templates for these flows. In particular, we describe how Ghys' modular template lifts to each  finite cyclic self-cover and determine the resulting knot types of the lifted periodic orbits. This gives a more geometric and visual interpretation of the self-covers studied in \cite{clay:20}.
\medskip

The trefoil has a Seifert surface $\Sigma$ homeomorphic to a once-punctured torus. By cutting the trefoil complement along this surface and considering the product $\Sigma \times [0,1]$, we construct a fundamental block of this space, as can be seen in figure 1.1 (b). The full complement is then recovered by gluing the top and bottom faces of this block using a specific monodromy - a rotation by $\frac{2\pi}{6}$ that we describe in detail in $\S$2. By stacking $n$ such blocks and applying the same rotational gluing between the top and bottom faces as described in \cite{rolfsen:03}, we obtain an $n$-sheeted cyclic self-cover of the trefoil complement, corresponding to one of the Anosov flows mentioned above. It follows that

\theoremm \emph{The embedding of the modular template in the trefoil complement is as depicted in Figure 1.1 (c). The template corresponding to an Anosov flow, which is the pre-image of the modular flow under the map $p_n$, is obtained by embedding one copy of the modular template in each of the $n$ blocks of the cyclic self-cover.}

\begin{center}
\begin{overpic}[scale=0.12]{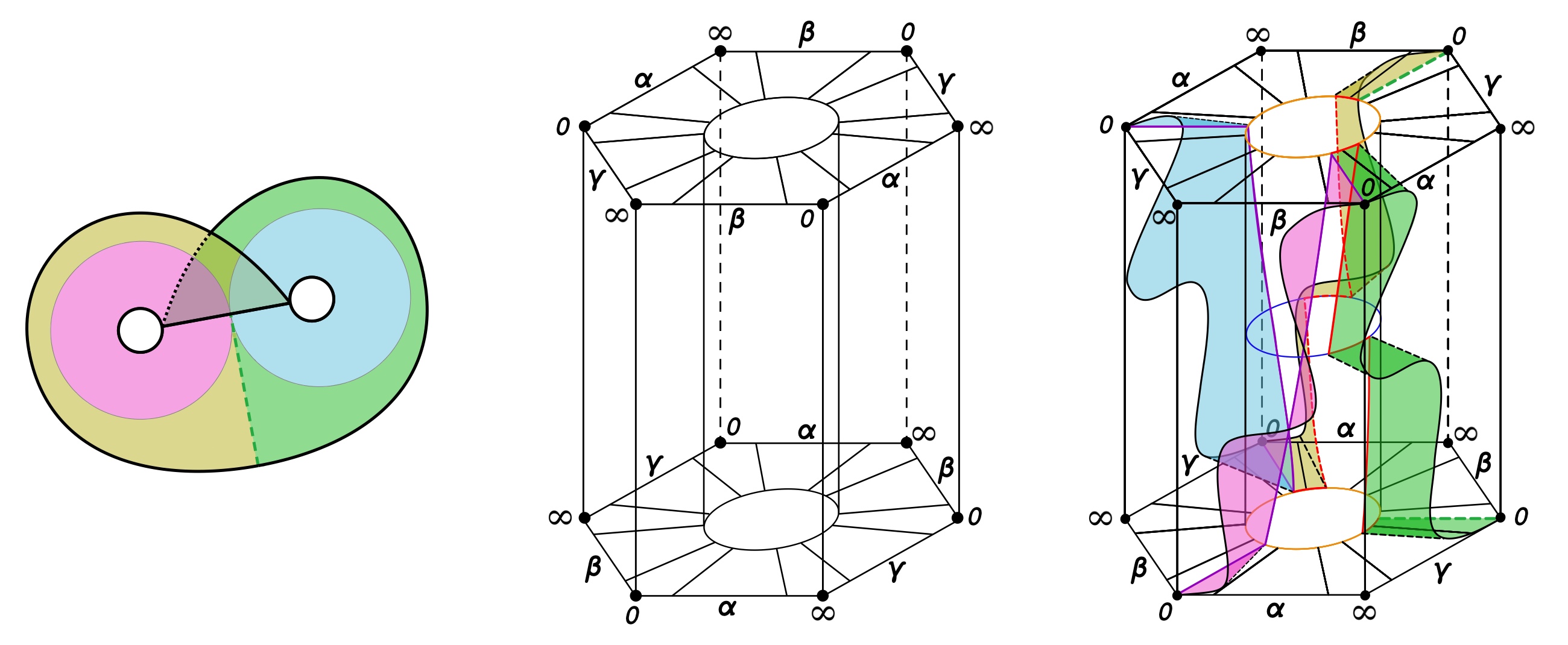}
\put(12.5,-3.5){\color{black}(a)}
\put(45,-3.5){\color{black}(b)}
\put(80,-3.5){\color{black}(c)}
\end{overpic}
\bigskip

\captionof{figure}{(a) The modular template; (b) A fundamental block representing the trefoil complement; (c) The modular template embedded in the trefoil complement}
\end{center}

The central result of our work is a complete description of how periodic orbits of the modular flow lift under these cyclic self-covers. In particular, we show that the lift of a periodic orbit represented by a word in the symbols $R$ and $L$ (introduced in $\S$2.2) decomposes into several connected components whose number depends explicitly on the exponents of $R$ and $L$. This statement, presented as our main theorem (Theorem 4.1.1), provides a precise $R,L$-coding of the lifted knots and establishes a clear combinatorial correspondence between the symbolic dynamics of the modular flow and the topology of the cyclic self-covers.
\medskip

This work also has a natural connection to commensurability. Friedl \cite{com:12} and Boileau \cite{boil:11} have both researched the commensurability of hyperbolic knots. Since we consider periodic orbits of the modular flow, which are hyperbolic knots in the trefoil complement \cite{foulon2013contact}, we conclude the following

\theoremm \emph{There are infinitely many periodic orbits $K$ on the modular template such that the link $K\cup T$ in $S^3$ is commensurable to infinitely many other links with two components, where one of them is $T$.}
\medskip

Finally, we discuss the implications of our construction for the arithmeticity of knot complements. As a consequence of the commensurability relation established in Theorem 1.2, we prove the existence of an infinite family of arithmetic knots in the trefoil complement, obtained by lifting modular knots that are themselves arithmetic (specifically, the $RL$ geodesic). This result is in contrast to the situation in $S^3$, where the figure-eight knot is the unique arithmetic knot.
\medskip

We now give a brief review of the overall scheme of the paper. In $\S$2, we provide information about this paper’s leading players: the trefoil knot and the modular template. We start by introducing the trefoil complement, its Seifert fibration, and their connected properties. Next, we describe the modular surface and its role in the geodesic flow and in creating the modular template. Then, we discuss finite cyclic self-covers of the trefoil complement. At the end of $\S$2, we build our model of the trefoil complement. In $\S$3, we begin sharing our research results. This section is fully dedicated to the embedding of the modular template in the model we build for the trefoil complement. The importance of this construction lies in understanding modular knots and in discussing their lifts to self-covers of the trefoil complement. After embedding the modular template, we state and prove the main theorem, Theorem 4.1.1, in $\S$4. As a final conclusion to this paper, we connect our results to commensurability and arithmeticity and prove Theorem 1.2. In $\S$5, we discuss open questions and future directions arising from our work.
\medskip

Throughout this paper, $n$ is always of the form $n=6k+1$ for $k\in \N_{>0}$, and every time we use the word rotation, it is in the counter-clockwise direction.
\medskip

We believe the work done in this paper would also help gain a clearer visual perception of the cyclic self-covers of the trefoil complement, and a better understanding of the behavior of periodic orbits of the geodesic flow.
\medskip

\textbf{Acknowledgments.} Above all, we wish to thank our supervisor, 
Prof. Tali Pinsky, for her guidance, patience and support throughout the research process. We would also like to thank the anonymous referee for their valuable comments and suggestions that have improved this paper. This work was inspired by the influential visualization done by Ghys and Leys, and we recommend watching the videos appearing in $\S3.1$ of \cite{Ghys:06}. 
\bigskip

\section{Preliminaries}

\subsection{Seifert fibration of the trefoil complement}
\medskip

The trefoil complement is a Seifert fiber space, that is, a decomposition of the trefoil complement in $S^3$ into infinitely many topological circles, which are one-dimensional fibers. The fibration includes infinitely many regular fibers, each isotopic to a trefoil in $S^3$, and two singular fibers, which are trivial knots. A nice visualization of the Seifert fibration of the trefoil complement appears in \cite{Ghys:06}.

\begin{center}
\begin{overpic}[scale=0.145]{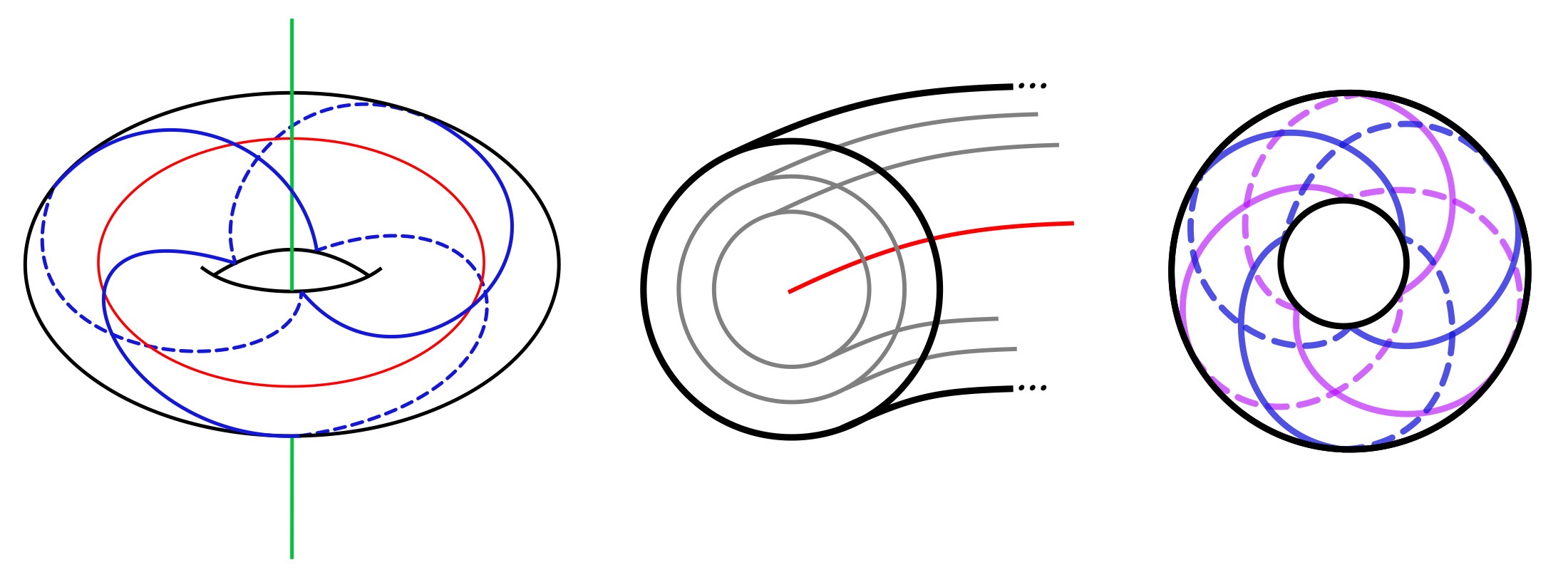}
\put(16.5,-3.5){\color{black}(a)}
\put(49,-3.5){\color{black}(b)}
\put(84.5,-3.5){\color{black}(c)}
\end{overpic}
\bigskip

\captionof{figure}{\text{Fibration of the trefoil complement: (a) A regular fiber (blue) and the two } \\ \text{singular fibers (red and green); (b) A cross section of the torus; (c) Two regular fibers }}
\end{center}

In addition to the one-dimensional fibration described above, another structure on the trefoil complement is a fibration over $S^1$, which arises from the action of $S^1$ on $S^3$ defined by the one-dimensional Seifert fibration. This structure, which we will refer to as the two-dimensional fibration of the trefoil complement, includes infinitely many fibers where each is a Seifert surface with a regular fiber as its one-dimensional boundary. We fix such a fiber and denote it by $\Sigma$, and mention that the Seifert surface of the trefoil is homeomorphic with the once-punctured torus. On the boundary of the trefoil complement, we define a basis consisting of a longitude and a meridian. The longitude is defined as the boundary of the Seifert surface $\Sigma$, while a meridian is a closed curve on the boundary of the trefoil that bounds a disk in the solid torus neighborhood of the trefoil and intersects the longitude exactly once. It is well-known that the trefoil complement $\cong\Sigma\times I$. Therefore, in the case of the two-dimensional fibration, we obtain a monodromy $\varphi: \Sigma \longrightarrow \Sigma$ by pushing forward every two-dimensional fiber (each point in $\Sigma$ is pushed forward along the one-dimensional fiber containing it) until it returns to $\Sigma$. We show in Lemma 2.3.1 that for an appropriate representation of $\Sigma$, $\varphi$ is a rotation by $\frac{2\pi}{6}$.
\medskip

In figure 2.1.2 (a) we observe a detailed projection of the Seifert surface of the trefoil. The yellow (resp. green) disk is a neighborhood of $0$ (resp. $\infty$), and we shall refer to it as the $0$-disk (resp. $\infty$-disk). The gray regions represent the twisted bands and the blue curve in both (a) and (b) is their boundary.
\medskip

Each fiber's linking number with the trefoil is the number of intersection points it has with the surface $\Sigma$. The regular fibers intersect the surface at six points, therefore each regular fiber has linking number $6$ with the trefoil. One of the singular fibers has linking number $3$ since it intersects the Seifert surface three times - once at each twisted band. The other singular fiber has linking number $2$ since it intersects the Seifert surface twice - once at the $0$-disk and once at the $\infty$-disk.
Let $f_t$ be the fiber with linking number $t$ with the trefoil. Notice that the longitude's linking number is $0$, as it is the boundary of $\Sigma$.

\begin{center}
\begin{overpic}[scale=0.24]{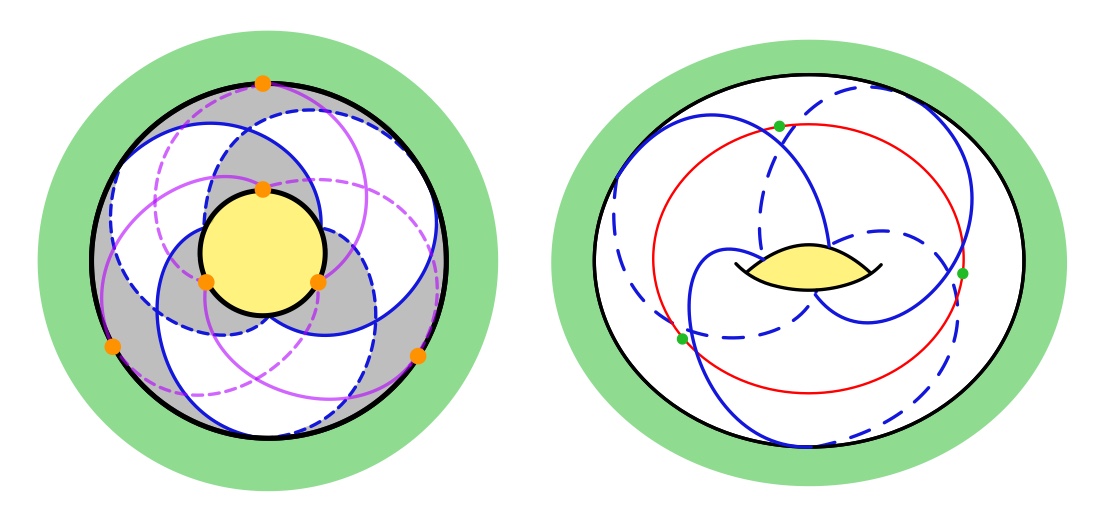}
\put(21.6,-3.5){\color{black}(a)}
\put(71.5,-3.5){\color{black}(b)}
\end{overpic}
\bigskip

\captionof{figure}{\text{(a) Intersection points of a regular fiber $f_6$ (orange);} \\ \text{(b) Intersection points of the singular fiber $f_3$ (green)}}
\end{center}

Later, we will consider the thickened trefoil on the torus and it will be necessary to understand how the parts of the Seifert surface $\Sigma$ lie on it. Thus, we give also the following figure, which exhibits the different parts with color correlation to figure 2.1.2 (a).

\begin{center}
\begin{overpic}[scale=0.15]{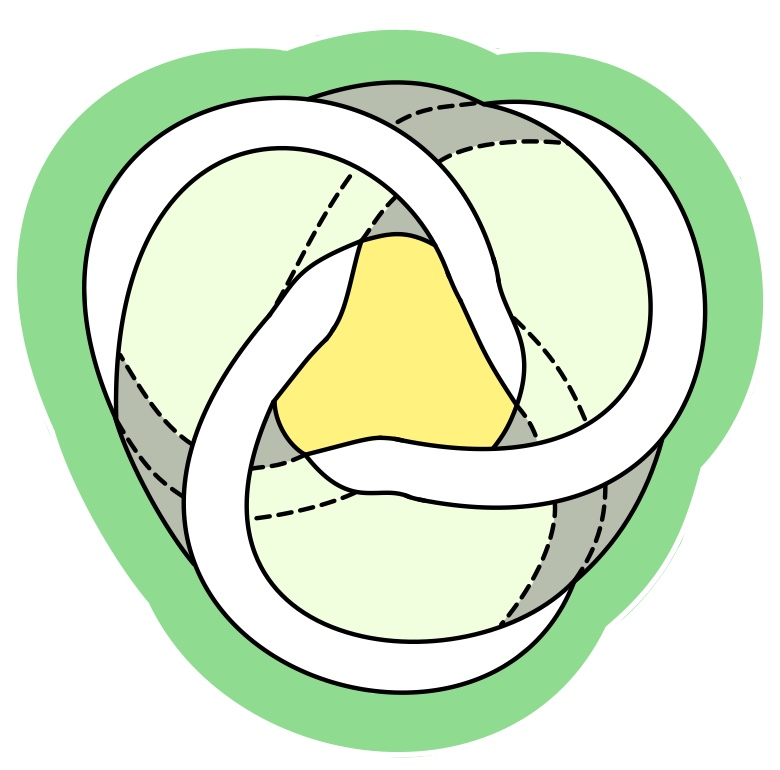}
\end{overpic}
\captionof{figure}{The Seifert surface $\Sigma$ on the torus}
\end{center}

Before moving on to the next subsection, we introduce an essential proposition on which we rely during our construction of the cyclic self-cover of the trefoil complement.
\medskip

\prop \emph{The monodromy $\varphi$ is of order $6$.}
\medskip

\proof  We only bring the important notions argued along this proof, which will be useful for us later. A full proof is given in \cite{rolfsen:03}, and the reader may see also \cite{clay:20}.
Recall that a fiber $f_6$ has linking number $6$ with the trefoil, and if we take two consecutive intersection points of an $f_6$ fiber with $\Sigma$, then they differ by a regular interval, which constitutes a sixth of the fiber. As a result, $\varphi$ cyclically permutes the intersection points along that fiber. Consequently, an intersection point with $f_6$ returns to itself only after acting $\varphi$ six times. This will result in permuting the $0$-disk and the $\infty$-disk, and in cyclically permuting the three twisted bands. Since $2$ and $3$ are coprimes, the intersection points of $f_2$ and $f_3$ with $\Sigma$ return to themselves only after acting $\varphi$, $\text{lcm}(2,3)=6$ times. In conclusion, for every point on $\Sigma$ to return to itself, we must act $\varphi$ exactly six times, hence $\varphi^6 =\text{Id}$. \hfill $\square$
\bigskip

\subsection{Constructing the modular template}
\medskip

In this subsection, we will construct the modular template from the point of view of a geodesic flow on the hyperbolic plane, and more specifically - on the modular surface, which will be introduced now.
\medskip

We identify the group of orientation-preserving isometries of the hyperbolic plane, denoted $\text{Isom}^+(\bb{H}^2)$, as $PSL_2(\R)$, acting on the upper-half plane by Möbius transformations. Let $\Gamma = PSL_2(\Z)$ be the modular group, a discrete subgroup of $\text{Isom}^+(\H^2)$, generated by

$$\Gamma = \left< s := \mat{0 & -1 \\ 1 & 0} ,  t := \mat{0&-1 \\ 1&-1 } \right>$$

The generator $s$ of $\Gamma$ acts as rotation by $\pi$ about the point $i$, making $i$ a cone point of order $2$ in the quotient orbifold $\H^2/\Gamma$, and the generator $t$ acts as a rotation by $\frac{2\pi}{3}$ about the point $q :=\frac{1}{2} + \frac{\sqrt{3}}{2}i$, making it a cone point of order $3$ in $\H^2/\Gamma$.
\medskip

The action of $\Gamma$ on ideal triangles in $\bb{H}^2$ gives rise to the Farey tessellation $\mathcal{F}$ of the hyperbolic plane. Correspondingly, we obtain a fundamental domain $\mathcal{D}$ for this action.

\begin{center}
\begin{overpic}[scale=0.3]{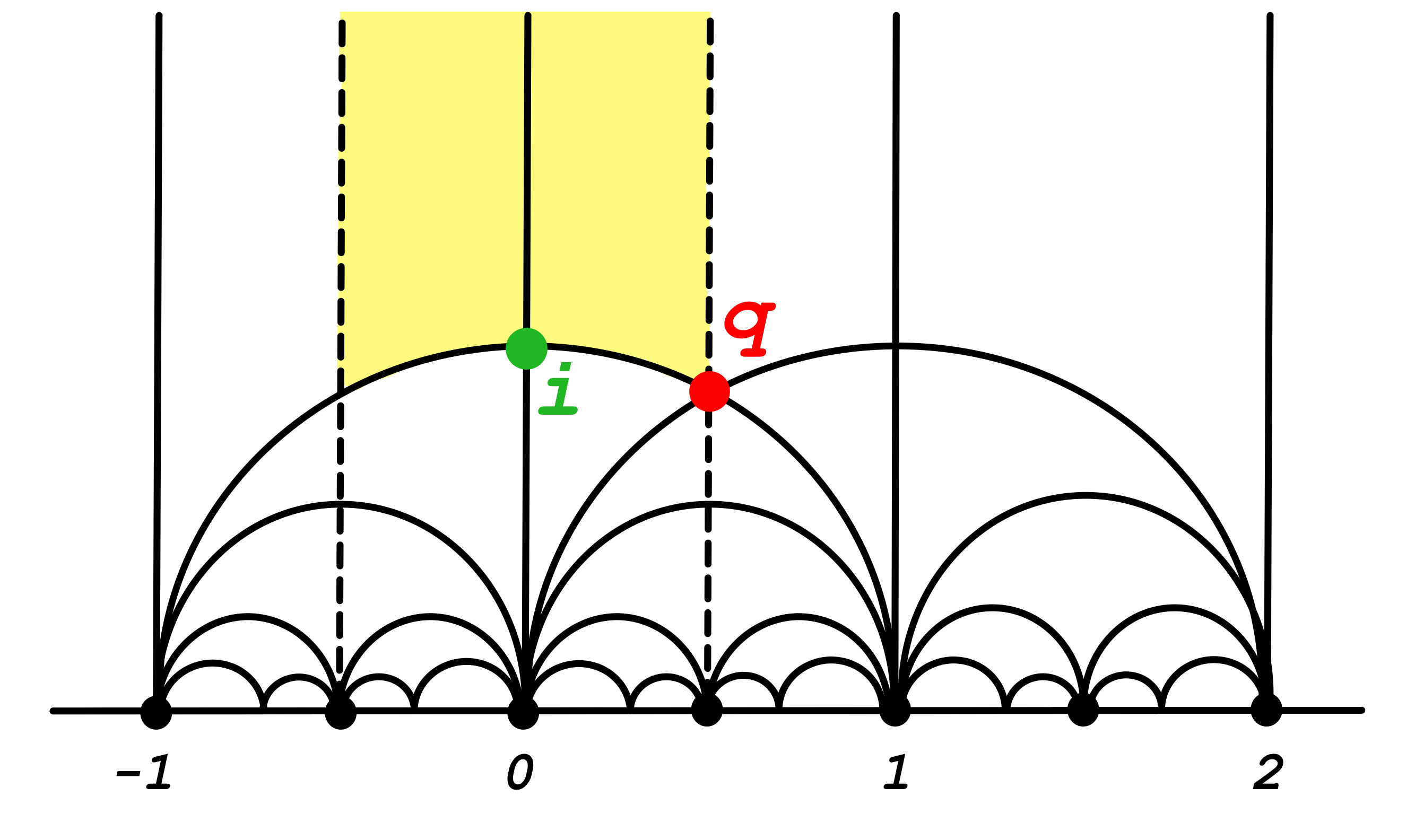}
\end{overpic}
\captionof{figure}{The Farey tessellation of $\H^2$, and the fundamental domain $\mathcal{D}$ (yellow)}
\end{center}
\medskip

The quotient space
$$S_{\text{Mod}} : = \H^2 /\Gamma$$
is \emph{the modular surface}. It is an orbifold with one cusp at infinity and two cone points - one of order $2$ and one of order $3$.

\begin{center}
\begin{overpic}[scale=0.04]{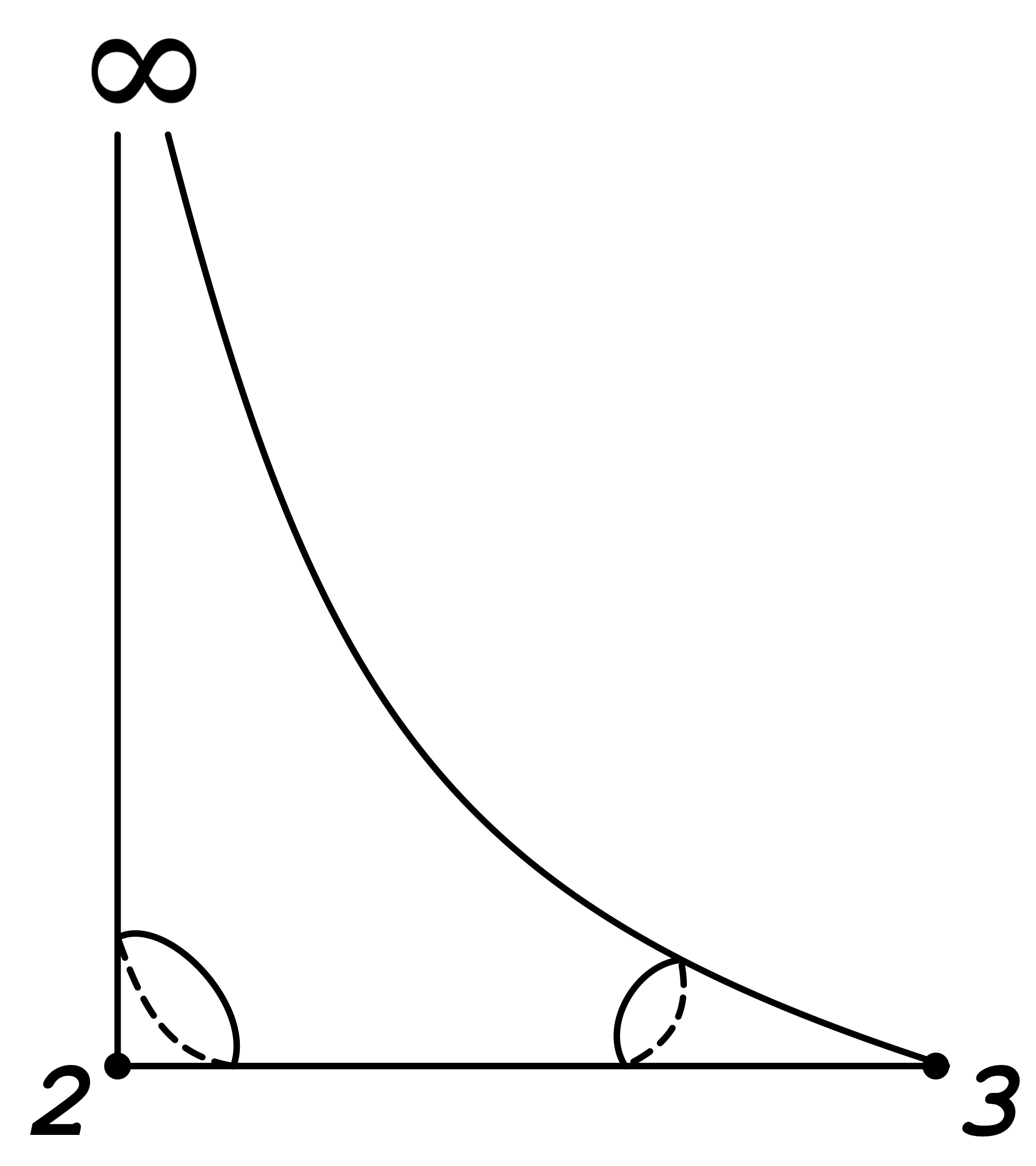}
\end{overpic}
\captionof{figure}{The modular surface $S_{\text{Mod}}$}
\end{center}

From the construction of $\mathcal{F}$, one observes that when a geodesic crosses an ideal triangle in $\mathcal{F}$ without terminating at an ideal vertex, it produces an oriented segment connecting two edges of that triangle. These two edges share a vertex;  we label the segment $R$ if the vertex lies to the right of the segment, and $L$ if the vertex lies to its left, as shown in figure 2.2.3. 
Since triangles in $\cF$ have vertices at $\Q \cup \brac{\infty}$, every geodesic in $\bb{H}^2$ with irrational endpoints is divided into segments and has an associated infinite itinerary. This itinerary is encoded by a sequence of the form $...R^{n_0}L^{n_1}R^{n_2}..., \,\, n_i \in \N$, which is invariant under $\Gamma$. For any geodesic $\gamma$ on $S_{\text{Mod}}$, we may lift $\gamma$ to a geodesic $\tilde{\gamma}$ in $\bb{H}^2$ and obtain a sequence of the above form. Different lifts of $\gamma$ differ by translations, which preserve orientation and leave $\cF$ invariant. Consequently, the label of each segment is preserved, and the sequence obtained is independent of the chosen lift. If  a geodesic on the modular surface is closed, there exists an element of $\Gamma$ translating points a finite distance along this geodesic, implying the geodesic has a periodic itinerary. See also \cite{Series:85} and \cite{Rickards:23}.
\medskip

\begin{center}
\begin{overpic}[scale=0.143]{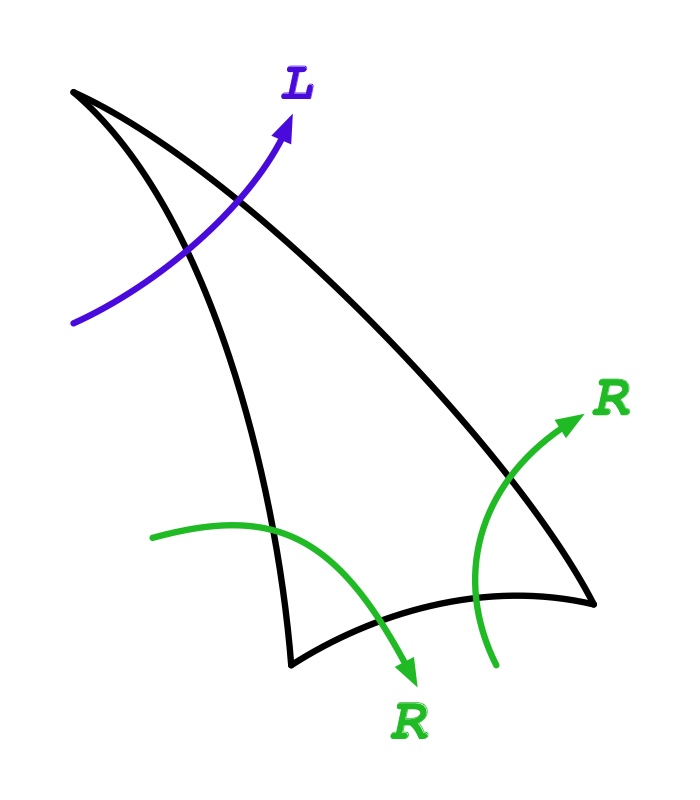}
\end{overpic}
\captionof{figure}{Labeling of segments as $R$ and $L$}
\end{center}

The construction of the modular template follows from the work of Ghys \cite{Ghys:07}, who established that the symbolic dynamics of geodesics in $\H^2$ provide a direct combinatorial representation of the flow on the branched surface. In this framework, the oriented segments crossing ideal triangles in the Farey tessellation correspond to the ears of the template according to their labels. By mapping these segments back into the fundamental domain $\mathcal{D}$ using  the isometries $L=st$ and $R=st^{-1}$ from $PSL_2(\Z)$, one ensures a coherent flow on the quotient surface. This allows us to explain how geodesics are projected onto $\mathcal{D}$ to form closed loops as follows.

As mentioned, every geodesic is composed of oriented segments crossing ideal triangles in $\mathcal{F}$. By acting on these triangles with isometries, we map each triangle, and consequently every contained segment, into $\mathcal{D}$. Up to isotopy, these segments do not intersect. This process ensures that all segments leave the same edge and are oriented outwards toward the right or left vertical edge of $\mathcal{D}$. Once completing this process, we obtain a fundamental domain containing oriented segments, where each segment carries a label $R$ or $L$. By gluing the marked edges as in figure 2.2.4, we obtain the modular template in its familiar form. Every periodic orbit corresponds to a Lorenz knot passing through the right and left ears of the template according to the geodesic's original itinerary. Note that, counter-intuitively, passing through the left ear of the template corresponds to $R$, while passing through the right ear corresponds to $L$.
\medskip

\begin{center}
\begin{overpic}[scale=0.047]{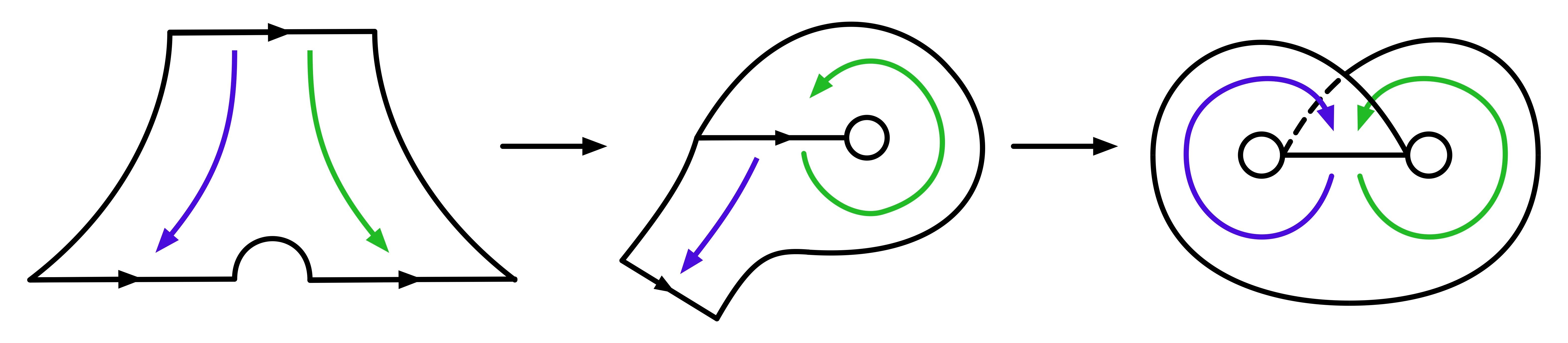}
\put(15.3,-3.5){\color{black}(a)}
\put(49,-3.5){\color{black}(b)}
\put(84,-3.5){\color{black}(c)}
\end{overpic}
\bigskip

\captionof{figure}{The modular template after gluing the edges}
\end{center}

The line along which the edges are glued is called the \emph{branchline}, as the template created is a branched surface.
\medskip

\conc \emph{Every periodic orbit in $S_{\text{Mod}}$ corresponds to a set of segments in $\mathcal{D}$, hence to a closed geodesic on the modular template. This gives rise to the desired correspondence between the geodesic flow and periodic orbits on the modular template \cite{brandts:19}.}
\medskip

The simplest example of a closed geodesic on the modular template is the $RL$ geodesic, depicted in figure 2.2.5.

\begin{center}
\begin{overpic}[scale=0.03]{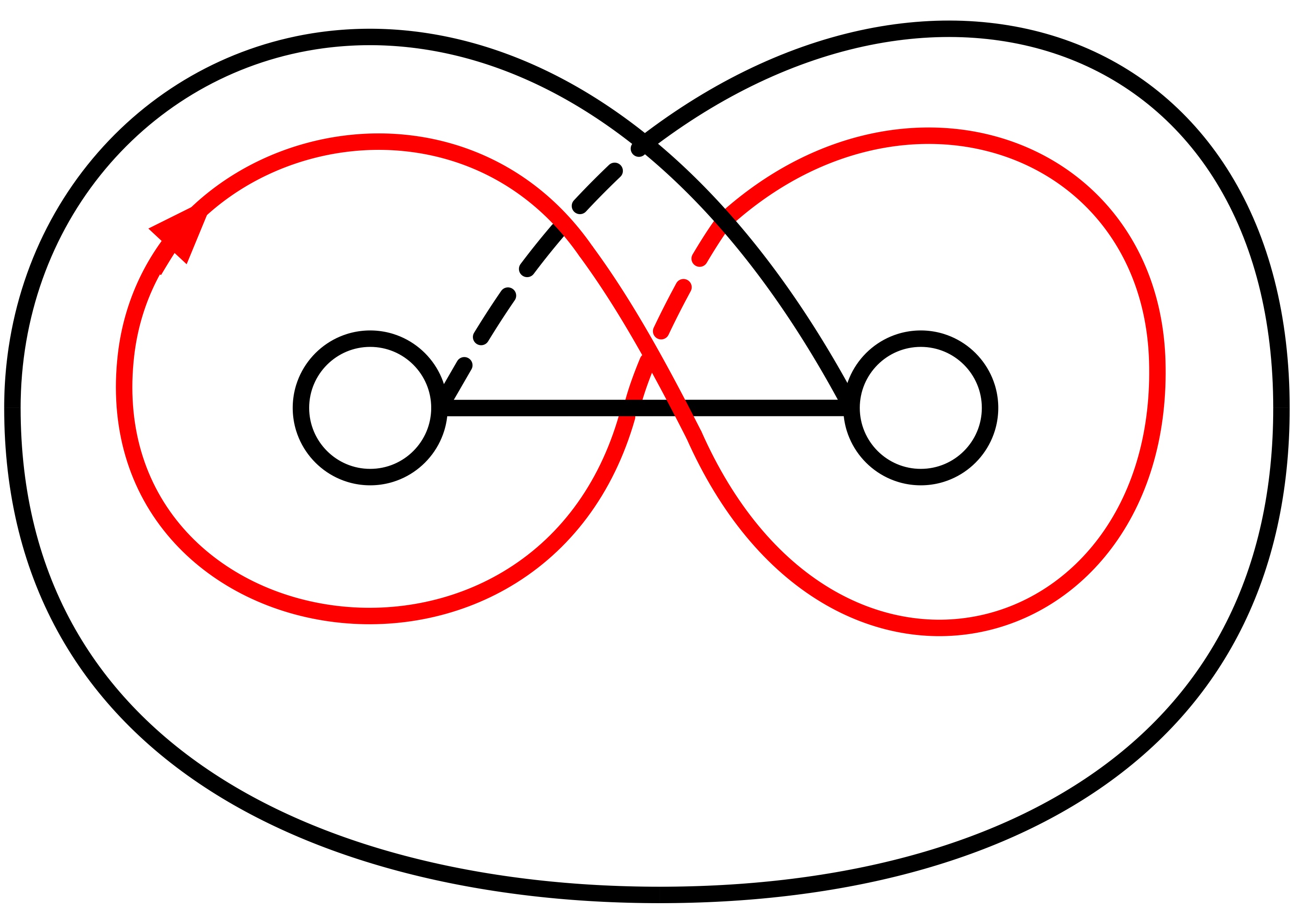}
\end{overpic}
\captionof{figure}{The $RL$ geodesic}
\end{center}

Ghys proved the equivalency between the symbolic dynamics of the representation of a geodesic on the modular template and its representation as a word in the generators $s,t \in PSL_2(\Z)$. Another example appears in figure 2.2.6.
\medskip

\begin{center}
\begin{overpic}[scale=0.03]{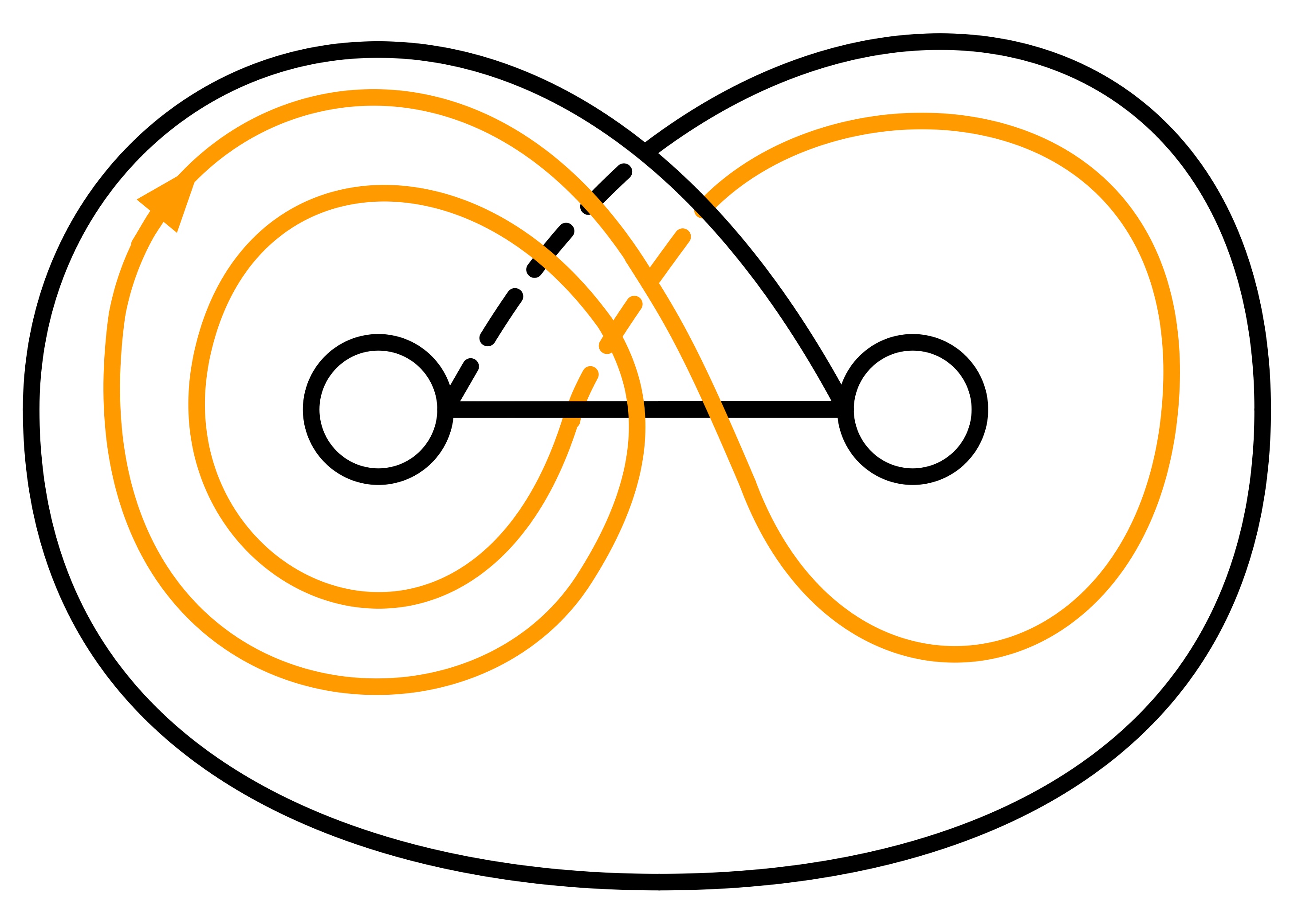}
\end{overpic}
\captionof{figure}{The $R^2L$ geodesic}
\end{center}

In conclusion, a natural bijection exists between closed orbits on the modular template and conjugacy classes of hyperbolic elements in $PSL_2(\Z)$. The primary purpose of our research is to construct self-covers of the trefoil complement and to discover which periodic orbits on the modular template lift to closed curves in these self-covers. Our technique will be to embed the modular template, its boundary, and the closed orbits corresponding to words in the symbols $R$ and $L$ in a self-cover of the trefoil complement, which we build in the following subsection. As a result, we will be able to embed every closed orbit in our model by following its corresponding word.

\subsection{Constructing a cyclic self-cover of the trefoil complement}
\medskip

We now explain our construction of the cyclic cover of the trefoil complement, according to \cite{rolfsen:03}. It follows from Proposition 2.1.1 that the cyclic cover of order $n$ is a self-cover.
\medskip

As mentioned, the Seifert surface homeomorphic with the once-punctured torus is the Seifert surface inducing the $S^1$ fibration. Our model will demonstrate precisely how the action of $S^1$  sweeps $\Sigma$ along the $3$-space, resulting in gluing with $\varphi$. In this paper, the once-punctured torus is represented by a punctured hexagon, denoted $H$. We choose a hexagon in order to take advantage of the symmetry of order $6$. Note that the equivalence classes consisting of two antipodal midpoints of edges are a periodic orbit of order $3$, while the equivalence classes of vertices form a periodic orbit of order $2$. Thus, we shall first examine the action of $\varphi$ on $H$.

\begin{center}
\begin{overpic}[scale=0.05]{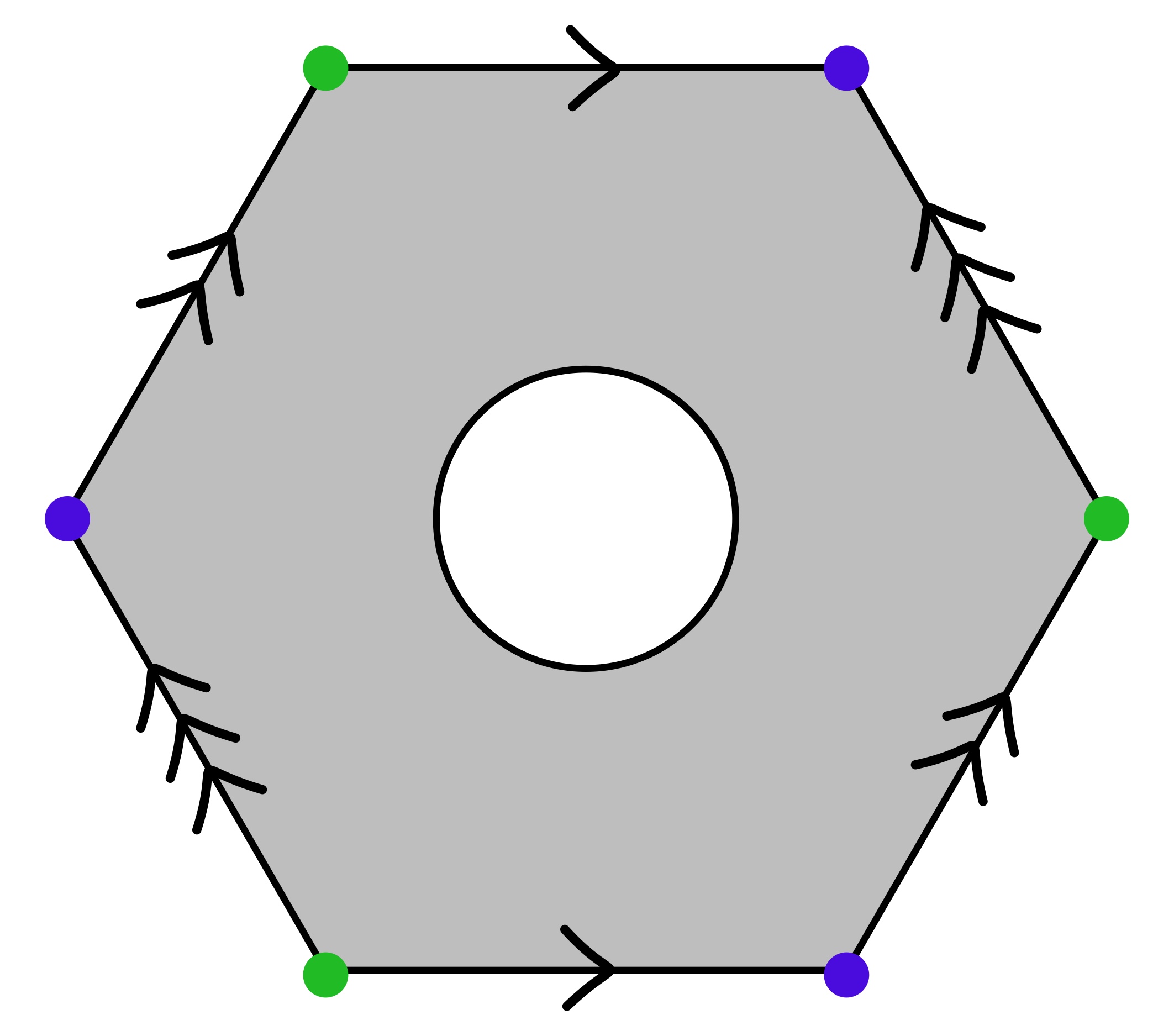}
\end{overpic}
\captionof{figure}{The once-punctured torus $H$}
\end{center}

\lemma \emph{The quotient of the once-punctured torus by the group generated by $\varphi$ is the modular surface $S_{\text{Mod}}$. In our chosen representation $H$, the action of $\varphi$ corresponds to a counter-clockwise rotation by $\frac{2\pi}{6}$.}
\medskip

\proof As established in Lemma 2.2 of \cite{pinsky-purcell2023}, the once-punctured torus forms a $6$-fold cover of the modular surface. The group of covering transformations is cyclic of order $6$, generated by a rotation of order three and a rotation of order two. To study this action, we choose a specific covering map and work explicitly with this model, representing the once-punctured torus as the hexagon $H$. In this model, the generator of the covering group acts as a rotation $\psi: H \to H$ by an angle of $\frac{2\pi}{6}$ counter-clockwise. Recall from the proof of Proposition 2.1.1 that the monodromy $\varphi$ exhibits the exact same permutation dynamics on the intersection points of the fibers: it permutes the $0$-disk and the $\infty$-disk, and cyclically permutes the three twisted bands. Because $\varphi$ generates this same cyclic action of order $6$, we can identify $\varphi$ with $\psi$. \hfill $\square$
\medskip 

\remark Note that $PSL_2(\Z)$ actually has two distinct subgroups of order $6$. While we use the one corresponding to a cover by the once-punctured torus, the other subgroup corresponds to a cover by the thrice-punctured sphere. These two distinct covers arise from the two different ways of gluing ideal triangles to form an ideal square in the Poincaré disk.
\medskip

We will now use the knowledge gained in the proof of Proposition 2.1.1 and Lemma 2.3.1 to understand the different parts assembling $H$. The $0$-disk's boundary is composed of three longitude components and three connections to each twisted band, as depicted in figure 2.3.2. We do the same construction for the $\infty$-disk as well.

\begin{center}
\begin{overpic}[scale=0.05]{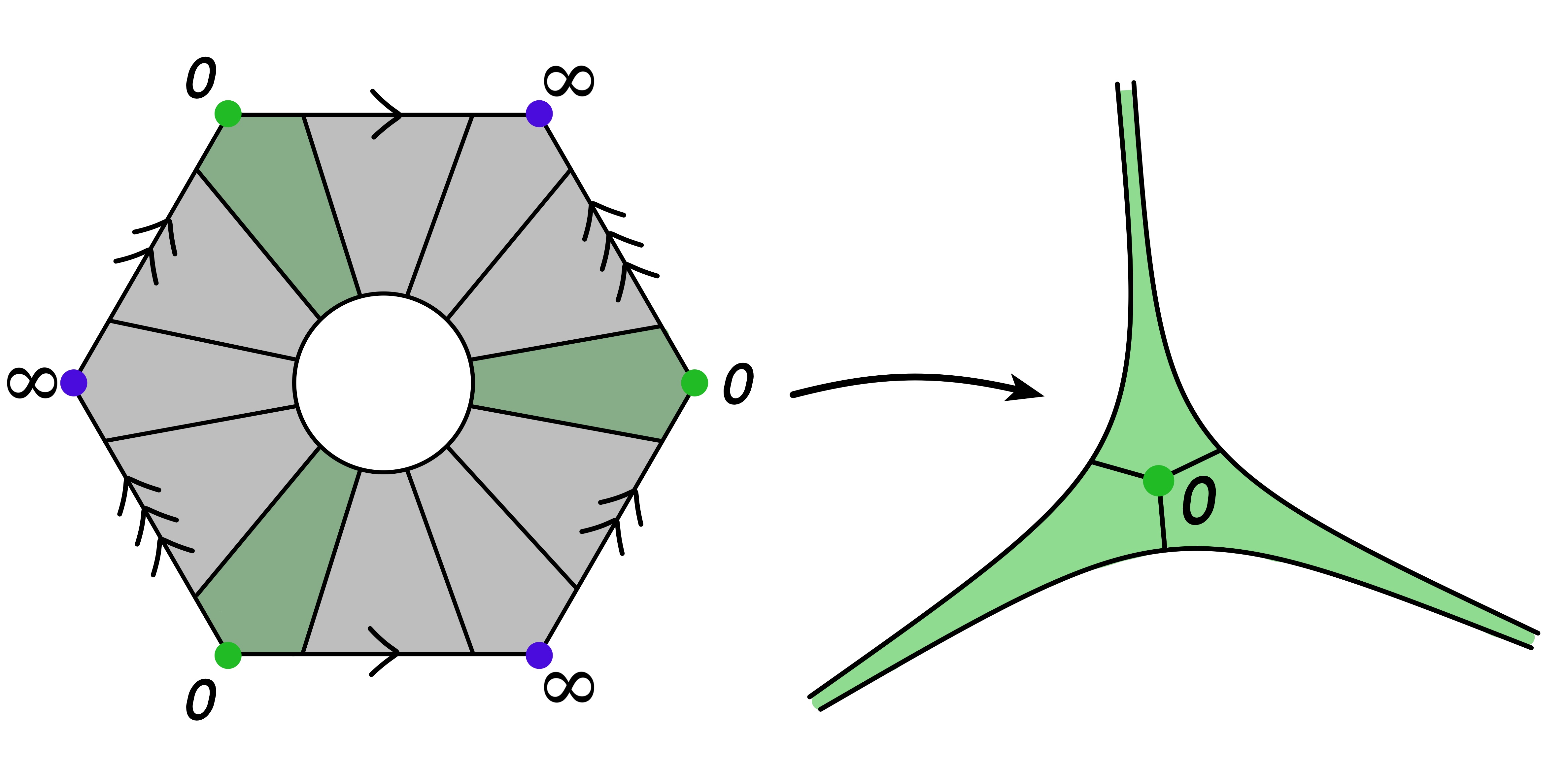}
\end{overpic}
\captionof{figure}{The $0$-disk, as a result of gluing the edges }
\end{center}

We also obtain the twisted bands $\alpha,\beta,\gamma$ by gluing every pair of opposite edges. In figure 2.3.3, every pair of regions adjacent to the mid-sections of edges of the same color represents one twisted band.

\begin{center}
\begin{overpic}[scale=0.05]{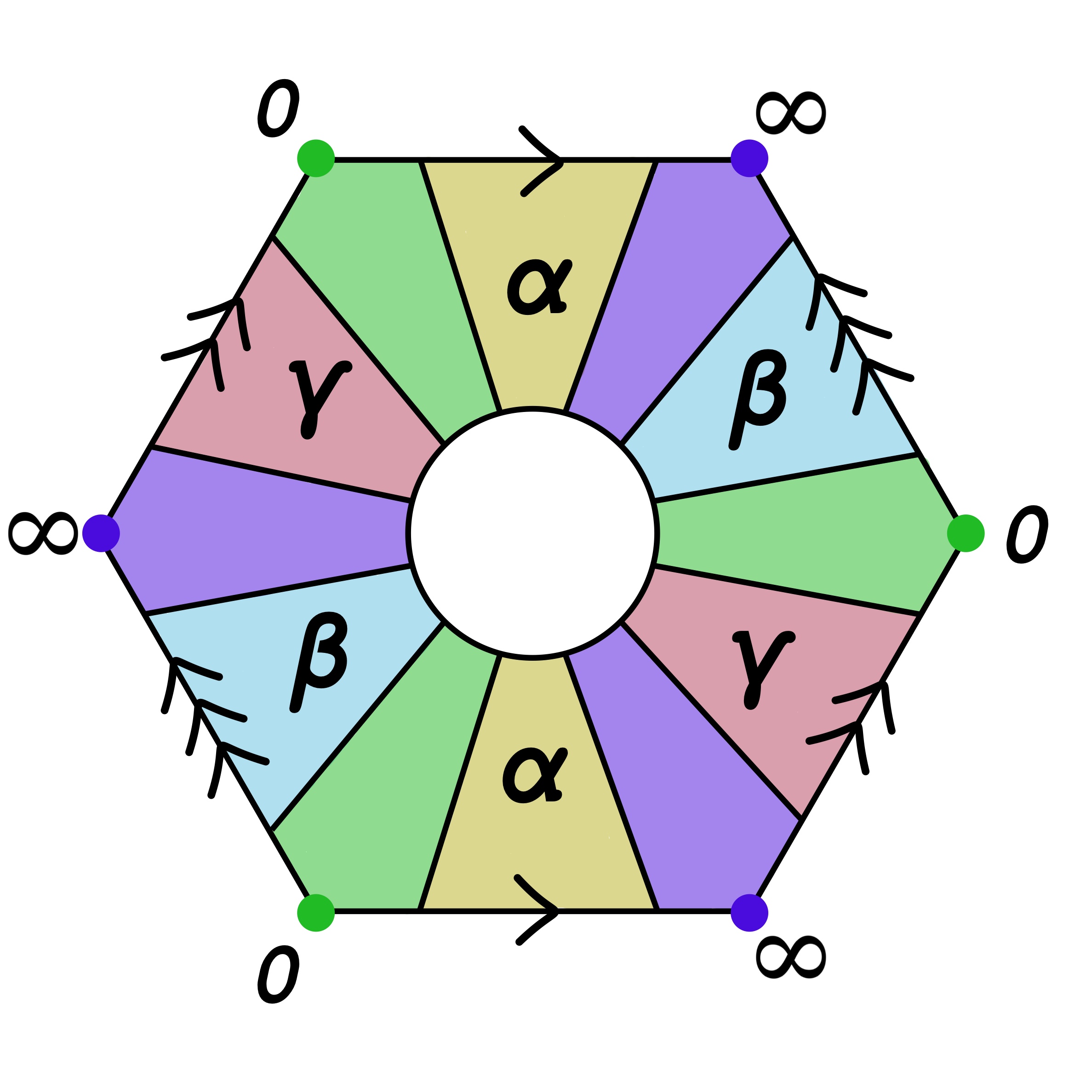}
\end{overpic}
\captionof{figure}{The twisted bands }
\end{center}

In figure 2.3.4, the blue edges represent the twisted band's boundary component, the red edge connects the twisted band with the $0$-disk, whereas the green edge connects the twisted band with the $\infty$-disk. We do the same construction for $\beta$ and $\gamma$ as well.

\begin{center}
\begin{overpic}[scale=0.05]{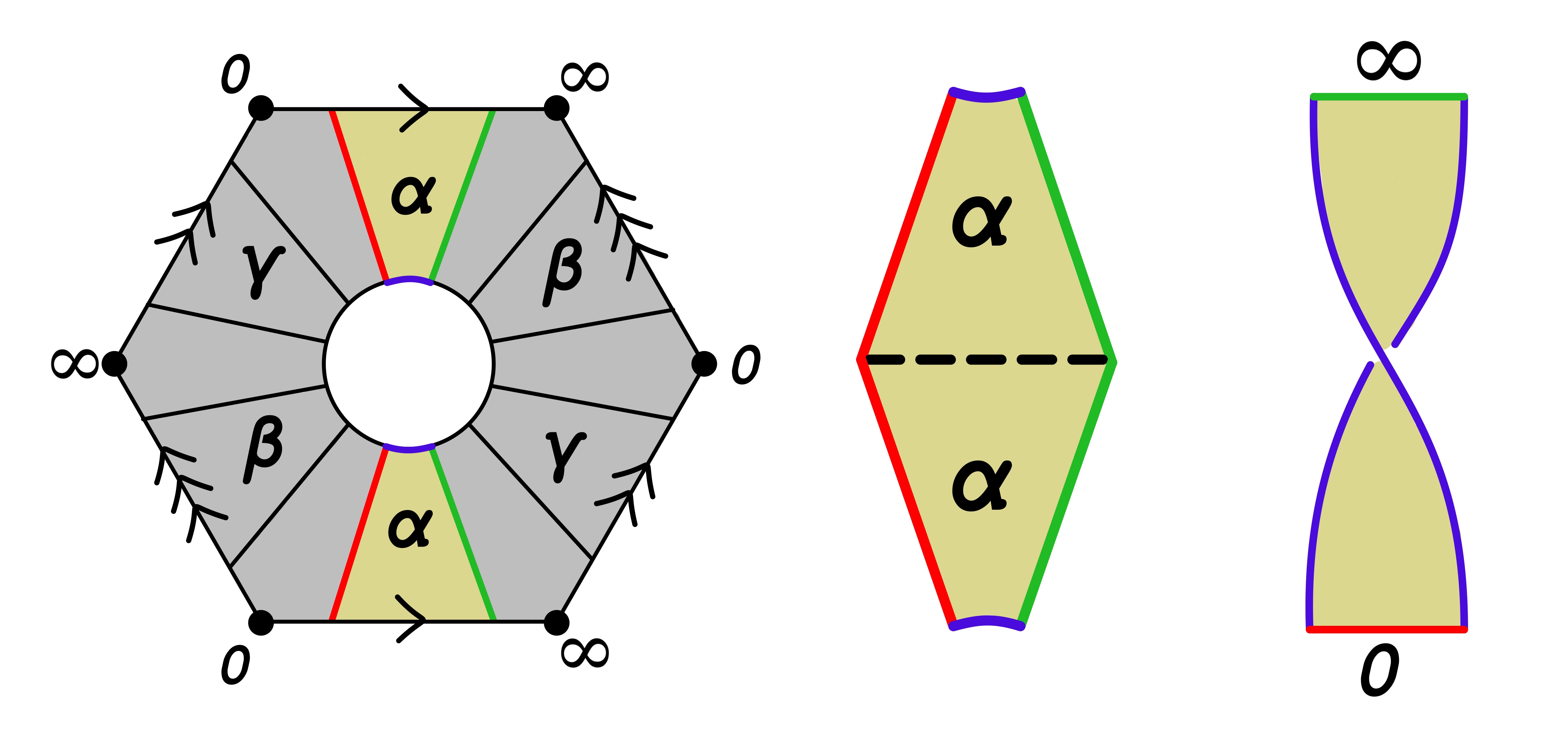}
\put(23.3,-3.5){\color{black}(a)}
\put(60.15,-3.5){\color{black}(b)}
\put(85.9,-3.5){\color{black}(c)}
\end{overpic}
\bigskip

\captionof{figure}{Twisted band $\alpha$; (a) before gluing; (b) after gluing; (c) after twisting}
\end{center}

Now, we are ready to look at the entire model.

\theorem \emph{The action of $S^1$ on the trefoil complement induced by \(\Sigma\) is described in the following figure}

\begin{center}
\begin{overpic}[scale=0.035]{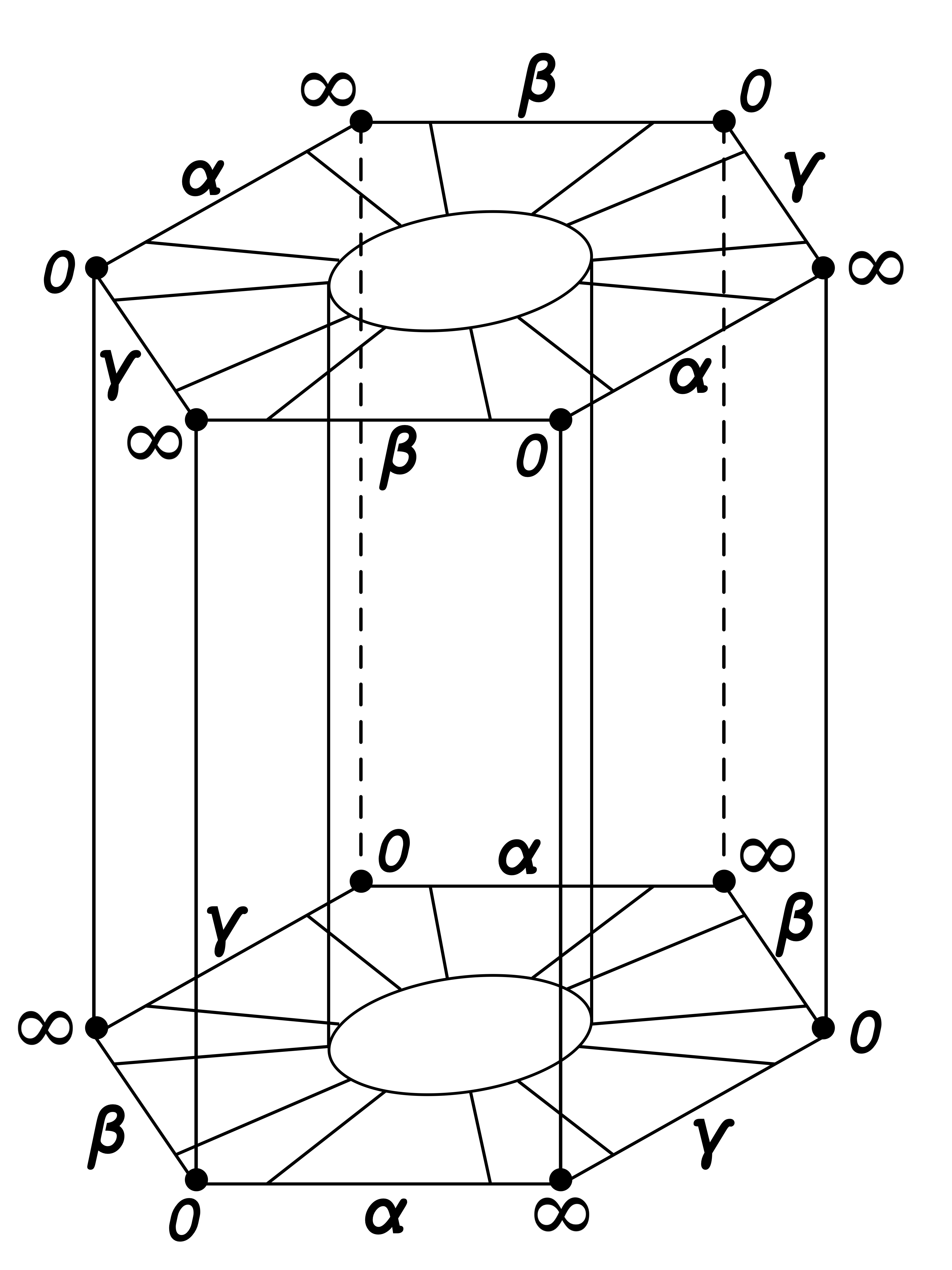}
\end{overpic}
\captionof{figure}{$\mathcal{H}_1$}
\end{center}

\proof The model $\mathcal{H}_1$ is constructed as  $(H \times I)/\sim$, where the equivalence relation identifies the top face $H\times \brac{1}$ with the bottom face $H\times \brac{0}$ via the monodromy $\varphi$. In this setting, the interval $I$ parametrizes the sweeping of the surface $\Sigma$ through the $3$-space under the action of $S^1$. Specifically, the interval $I$ corresponds to a single action of $\varphi$, which we have shown is equivalent to rotating the top face by $\frac{2\pi}{6}$. The identification of the top and bottom faces thus represents the return of the Seifert surface $\Sigma$ to its initial position under the rotational action. \hfill $\square$
\medskip

So far we understand the correspondence between each part of $\mathcal{H}_1$ to the sweeping of $\Sigma$ through the $3$-space under the action of $S^1$. Specifically, the longitude at $t\in I$ is the boundary of $H \times \brac{t}$. Thus, we can see clearly where the meridians lie within $\mathcal{H}_1$.
\medskip

\prop \emph{Every meridian appears in $\mathcal{H}_1$ as a vertical line along the middle cylinder, twisted by an angle of $\frac{2\pi}{6}$ counter-clockwise.}
\medskip

\begin{center}
\begin{overpic}[scale=0.035]{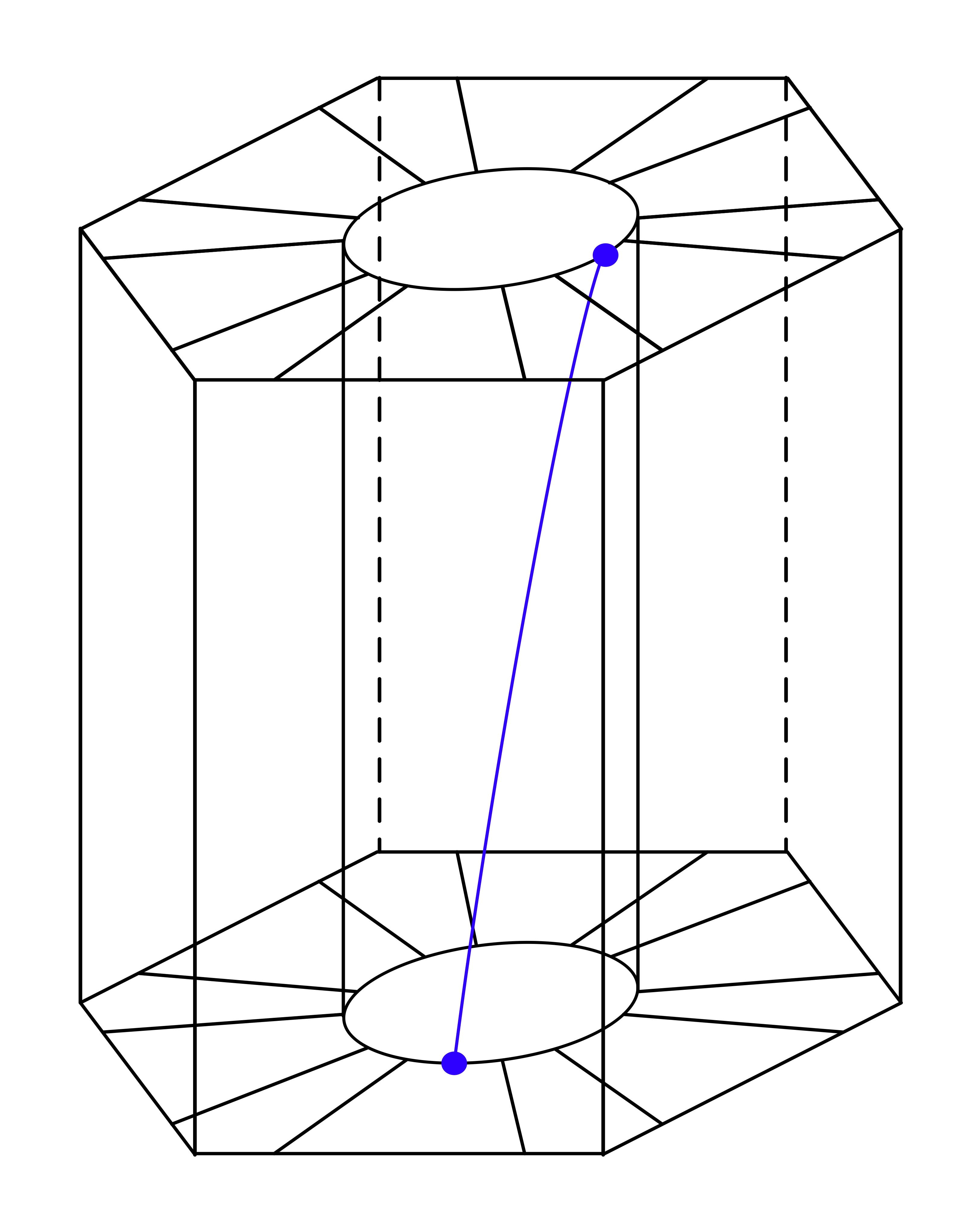}
\end{overpic}
\captionof{figure}{Meridian (blue)}
\end{center}

\proof Every meridian lies on the boundary and intersects a longitude exactly once, hence for each $t\in I$, the meridian is represented as a single point lying on the boundary of $H \times \brac{t}$. The top face is glued to the bottom face with a rotation of $\frac{2\pi}{6}$, therefore the two points representing the meridian in these faces are identified. Since the movement is continuous, the meridian is represented as a continuous line connecting the two edges. We purposefully choose the simplest such curve, omitting any unnecessary horizontal winding. \hfill $\square$
\medskip

\note \emph{We assign an orientation to the meridian and longitude as follows. The movement along the meridian is consistent with the action of $S^1$, i.e., while we go up the meridian, we twist with an angle of $\frac{2\pi}{6}$ counter-clockwise, and vice-versa. On the longitude, our path is in the clockwise direction: $\alpha \rightarrow \beta \rightarrow \gamma$ as can be seen in figure 2.3.7. The longitude and the meridian are essential to our work in the sense that this orientation allows us to use them as a coordinate system, through which we can describe the movement of the trefoil complement in the $\mathcal{H}_1$ model.}
\medskip

\begin{center}
\begin{overpic}[scale=0.035]{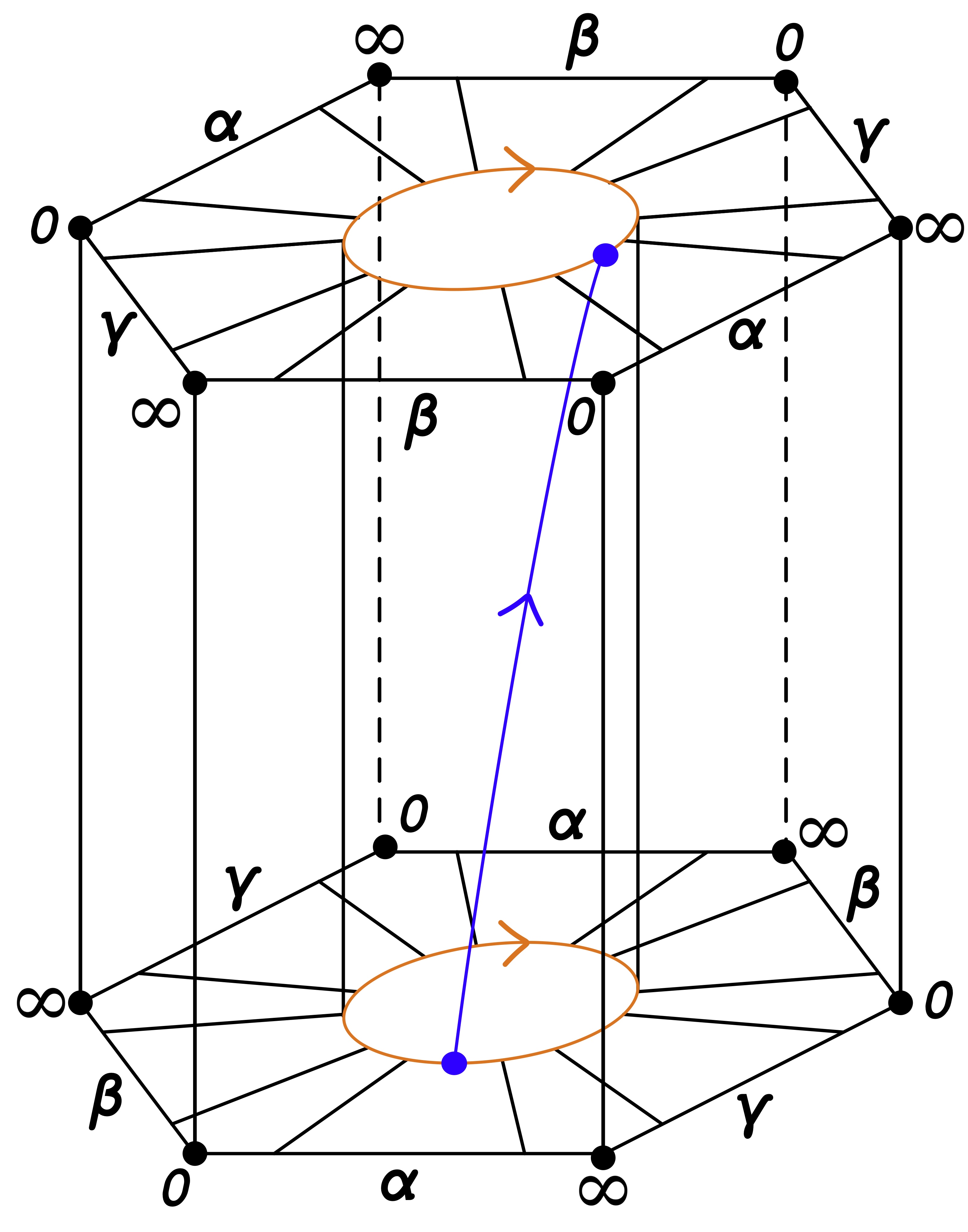}
\end{overpic}
\captionof{figure}{Meridian (blue) and longitude (orange), with orientation}
\end{center}

Every cyclic cover of order $n$ is constructed from $n$ pieces of $\mathcal{H}_1$. The bottom face of the following piece is glued vertically to the top face of the previous, such that all parts of these faces are aligned. After gluing the $n$ pieces together, the top face of the uppermost block and the bottom face of the lowermost block are also glued to one another with $\varphi$. 
\medskip

Now that we have a model of the trefoil complement sweeping along the fibered space and a coordinate system describing it, we can embed the modular template within the model.
\bigskip

\section{Embedding the template}

With the full terminology of the cyclic self-covers and the monodromy action established, we can refine the statement of Theorem 1.1 into two precise results: the explicit embedding of the template in a single block (Theorem 3.1) and its extension to the $n$-sheeted cover (Theorem 3.2).
\medskip

\theoremm \emph{The embedding of the modular template in $\mathcal{H}_1$ is as appears in figure 3.1.}

\begin{center}
\begin{overpic}[scale=0.2]{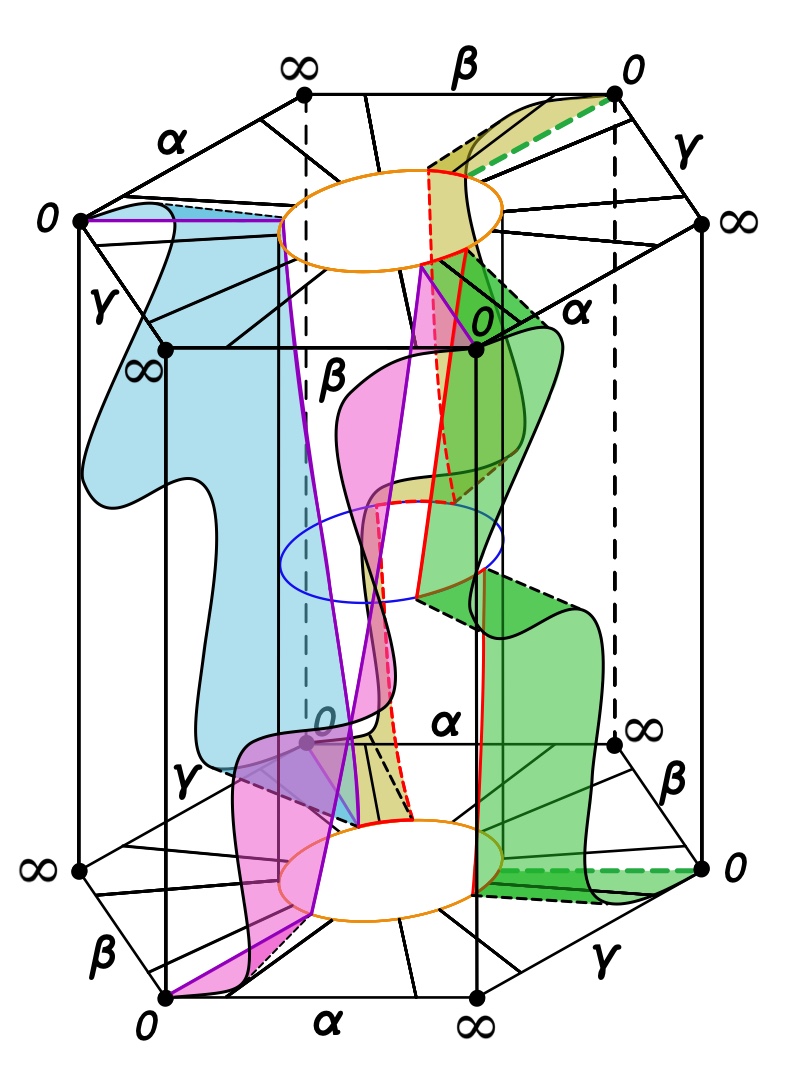}
\end{overpic}
\captionof{figure}{The modular template in $\mathcal{H}_1$}
\end{center}

\proof We start from the embedding of the modular template in the trefoil complement as computed by Ghys \cite{Ghys:07}. In particular, he showed that the inner boundaries of the template are meridians. This embedding together with the $RL$ geodesic can be seen in figure 3.2. Considering the $RL$ geodesic in this context would be beneficial in $\S$4.

\begin{center}
\begin{overpic}[scale=0.055]{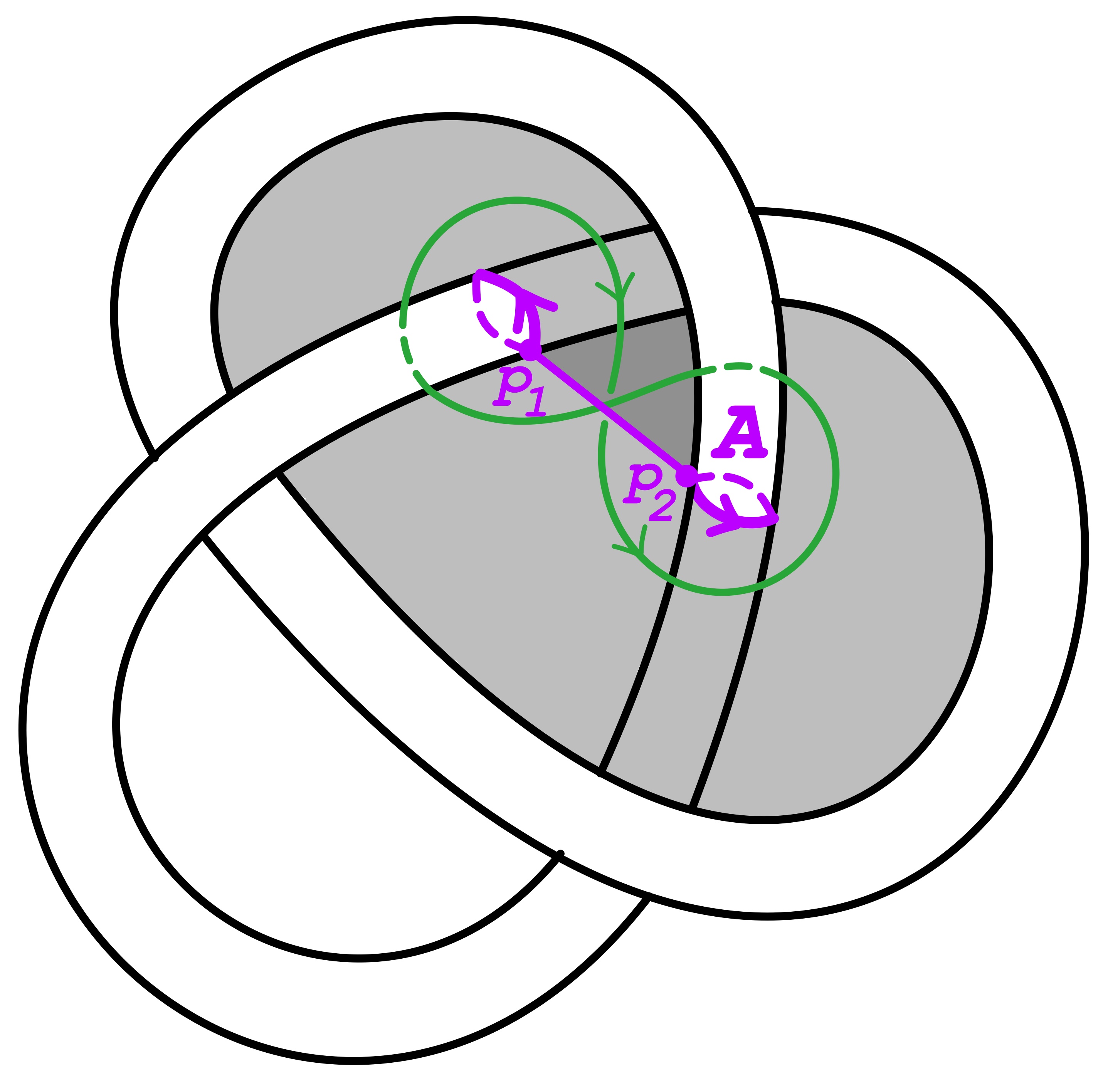}
\end{overpic}
\captionof{figure}{The modular template embedded in the trefoil complement}
\end{center}

As in figure 3.2, we denote by $A$ the purple set composed of the two meridians and the branchline that lies on the template's boundary. The green loop in this figure is the $RL$ geodesic as it was presented in $\S 2.2$, and we can continuously deform it into the set $A$. This is done by tightening the $R$ and $L$ factors into a small neighborhood of the thickened trefoil, such that the two loops are isotopic. Continuing to tighten the loop onto the thickened trefoil results exactly in the set $A$, which is homotopic to the original loop. The purple path created begins at $p_1$, travels a single meridian (which corresponds to the $x$ factor), then crosses the branchline in one direction until it meets $p_2$. Afterward, it goes through a single meridian (which corresponds to the $y$ factor) and crosses back to $p_1$ in the opposite direction of the branchline.
\medskip

We begin by embedding the template's boundary in $\mathcal{H}_1$, composed of the set $A$ (with one meridian in reversed orientation) and the remaining red part, which we denote by $B$. Both sets are depicted in figure 3.3.

\begin{center}
\begin{overpic}[scale=0.045]{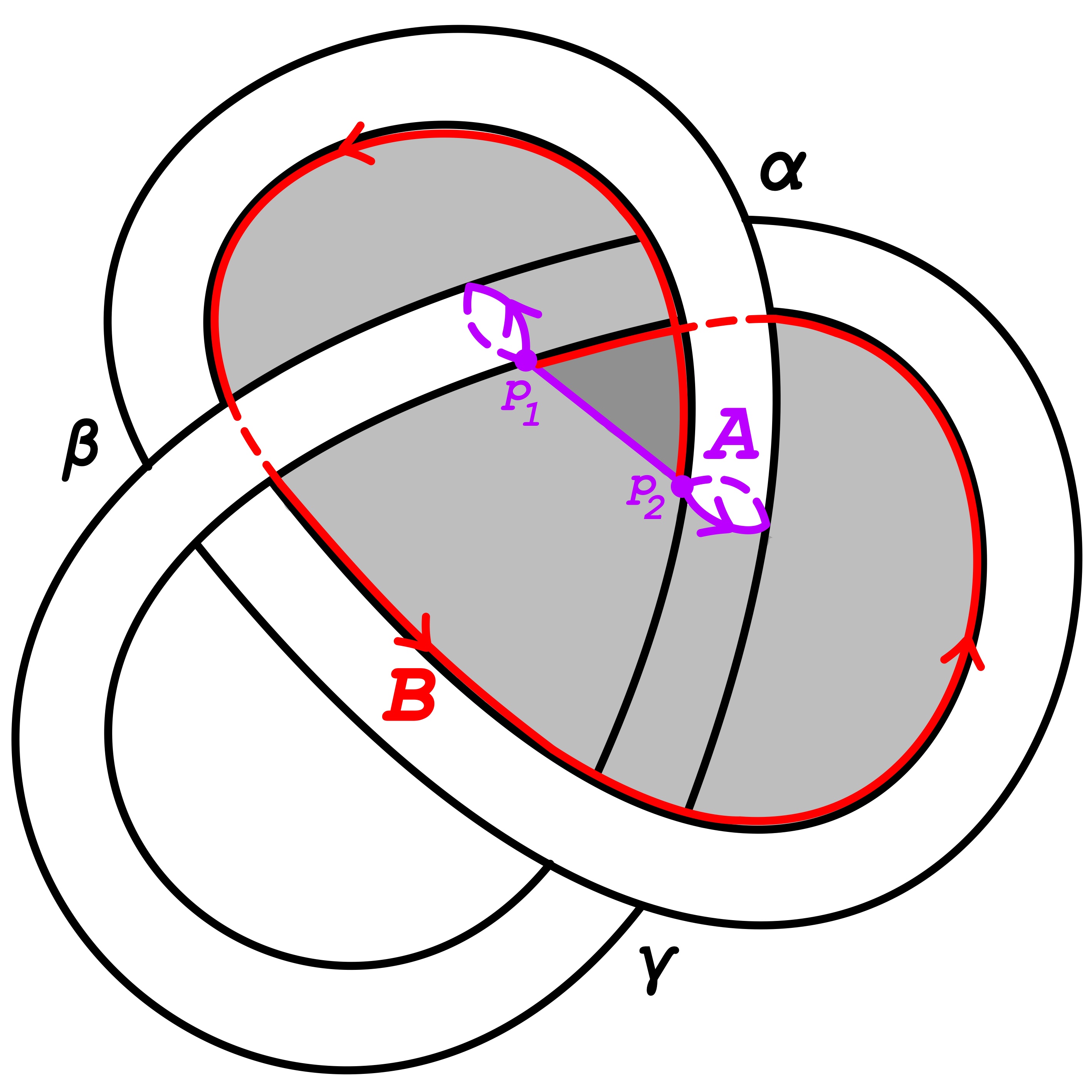}
\end{overpic}
\captionof{figure}{The template's boundary, composed of the sets $A$ (purple), and $B$ (red)}
\end{center}

The fact that $A$ is composed of meridians and the branchline allows us to prove its embedding in $\mathcal{H}_1$ is as in figure 3.4.

\begin{center}
\begin{overpic}[scale=0.035]{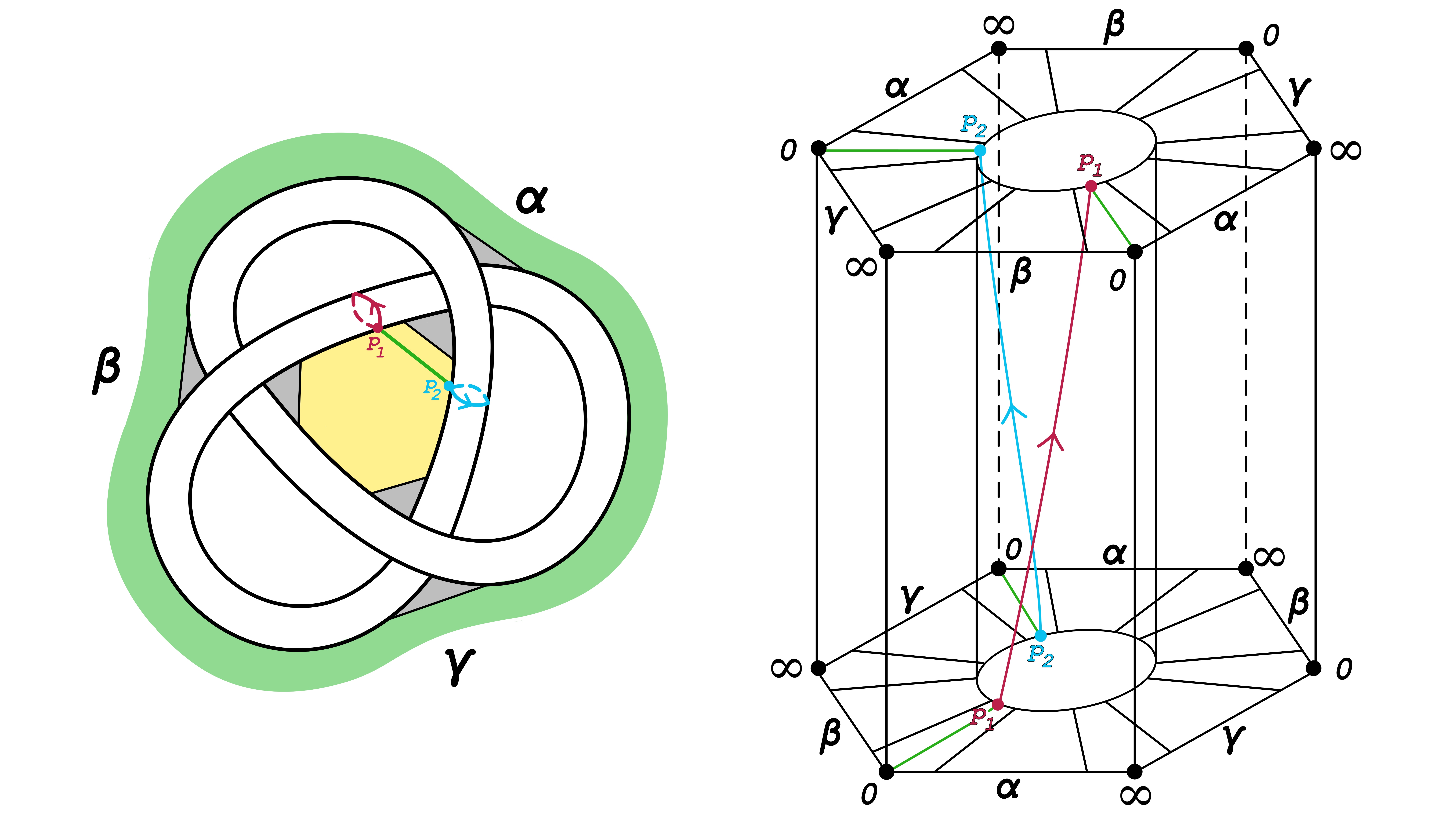}
\put(23.75,-3.5){\color{black}(a)}
\put(70,-3.5){\color{black}(b)}
\end{overpic}
\bigskip

\captionof{figure}{The embedding of $A$ in $\cH_1$}
\end{center}

As mentioned in $\S 2.1$, the yellow region in figure 3.4 (a) is the $0$-disk of $\Sigma$, so we may assume for convenience that the point $0$ is the midpoint of the branchline. Notice that $p_1$ is located between twisted bands $\alpha$ and $\beta$, and lies on the boundary of the $0$-disk where it connects to the thickened trefoil. Therefore, in $\cH_1$ it will appear in two copies as in figure 3.4 (b) (an identical argument holds for $p_2$). As for the meridians, followed by Proposition 2.3.4, each meridian is based at a point $p_i$ and connects the two copies. Note that in figure 3.4 we have already marked the orientation on the meridians according to Note 2.3.5 in both $(a)$ and $(b)$. Considering again figure 2.3.2, we can see that the branchline is composed of two straight segments, one that connects $0$ and $p_1$ and another one that connects $0$ and $p_2$. Both of them compose a straight line that is located exactly between twisted bands $\alpha$ and $\beta$. Next, recall that in $\mathcal{H}_1$ the $0$ points are all identified, therefore we may continue the line from the region located between twisted bands $\alpha$ and $\gamma$, and then we reach $p_2$. We draw the branchline in both top and bottom faces.
\medskip

We are able to choose two parallel longitudes on the thickened trefoil such that their union contains the set $B$ entirely, up to the boundary parts that pass through the twisted bands. The specific choice we made for $\lambda_1$ and $\lambda_2$ can be seen in figure 3.5, together with the boundary of the template composed of the sets $A$ and $B$.

\begin{center}
\begin{overpic}[scale=0.065]{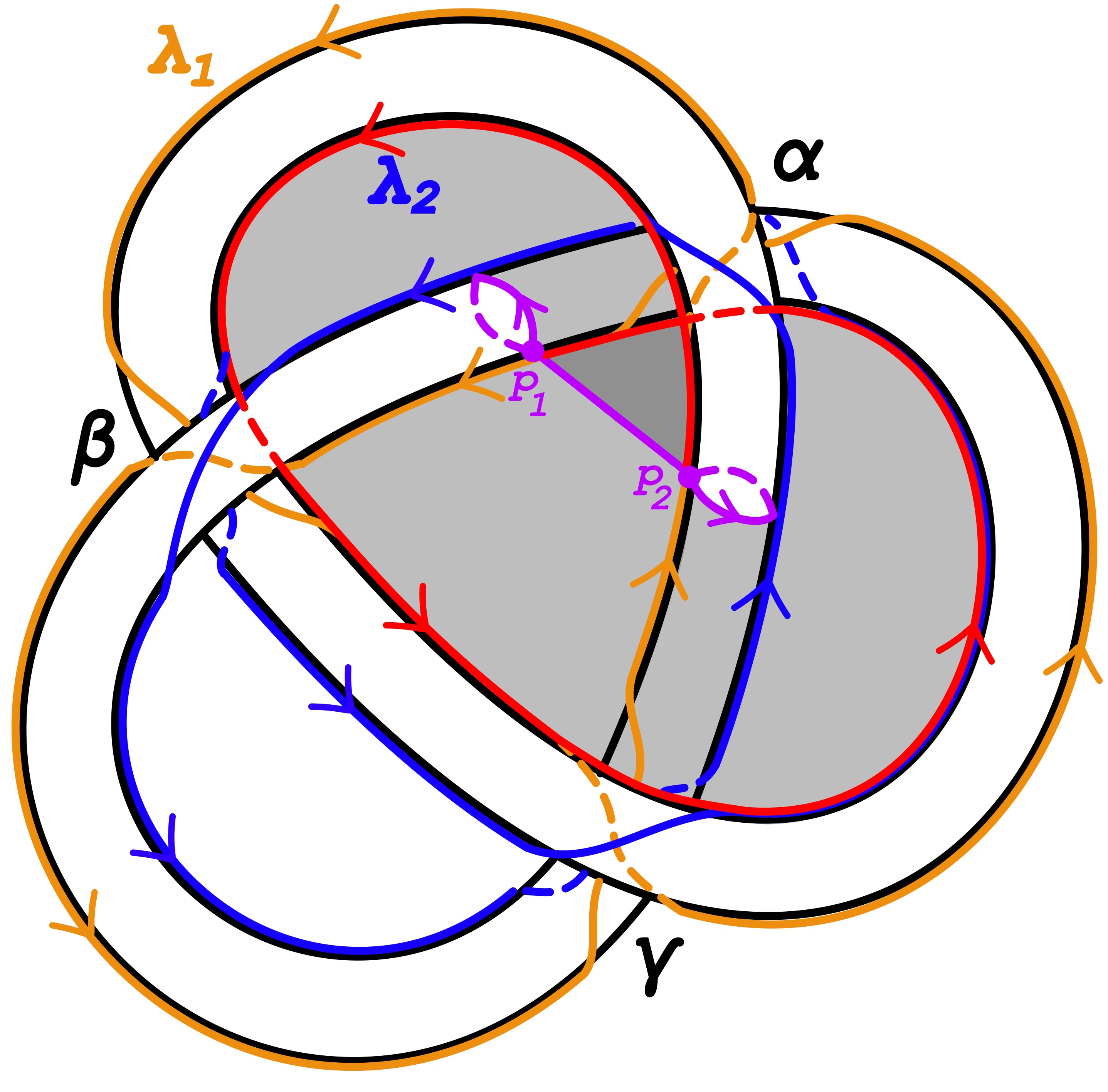}
\end{overpic}
\captionof{figure}{The modular template's boundary in the trefoil complement (red and purple), \\ and the two longitudes $\lambda_1$ (orange) and $\lambda_2$ (blue) containing the boundary}
\end{center}

We now prove that the embedding of the chosen longitudes, $\lambda_1$ and $\lambda_2$, in $\mathcal{H}_1$ is as depicted in figure 3.6.

\begin{center}
\begin{overpic}[scale=0.04]{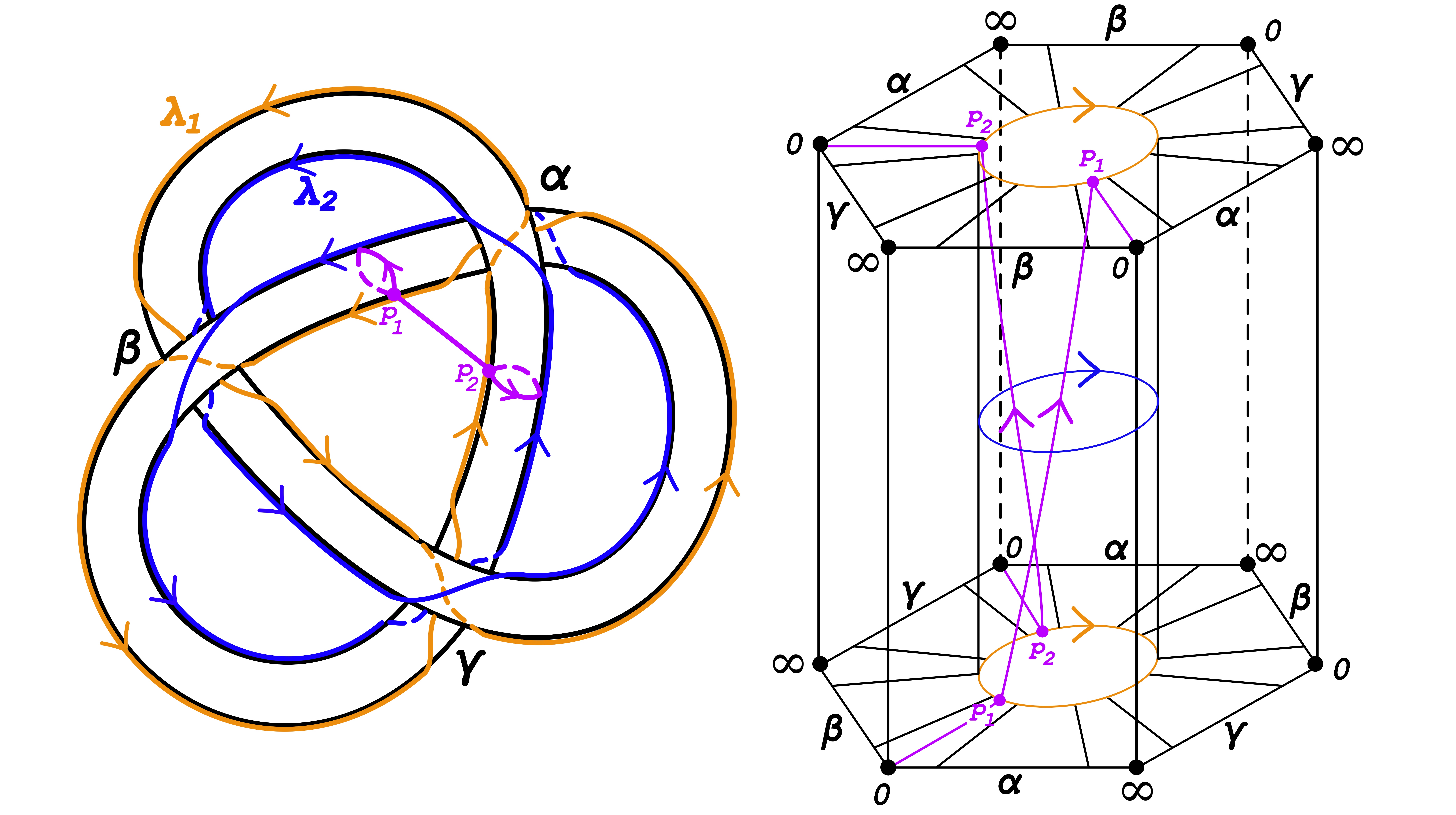}
\put(23.75,-3.5){\color{black}(a)}
\put(70,-3.5){\color{black}(b)}
\end{overpic}
\bigskip

\captionof{figure}{The longitudes $\lambda_1$ and $\lambda_2$, and the set $A$ in both depictions \\  of the trefoil complement }
\end{center}

First, the orientations on the $RL$ geodesic, on $\lambda_1$, and on $\lambda_2$ all correspond to the orientations as in Note 2.3.5.
As mentioned before, $\lambda_1$ is exactly the boundary of $\Sigma$; accordingly, it is represented by two orange circles on the top and bottom faces. $\lambda_1$ and $\lambda_2$ differ by an angle of $\pi$ along a meridian, so we reach $\lambda_2$ after half of the action of $\varphi$, hence it is represented by a blue circle in the mid-height of $\mathcal{H}_1$. Following figure 3.6 (a), we know exactly which parts of $\lambda_1,\lambda_2$ meet the modular template's boundary, as was shown in figure 3.5. 
\medskip

With that in mind, we can embed also the set $B$ in $\cH_1$. Notice that when $\lambda_1$ passes through twisted band $\alpha$ it goes \textbf{under} the thickened trefoil. At that point, $\lambda_2$ goes \textbf{over} it. Each transition is embodied in traveling minus half a meridian according to our chosen orientation. So if we set our starting point to be $p_2$, we travel along $\lambda_1$, and when we reach twisted band $\alpha$, we part from $\lambda_1$ and continuously coincide with $\lambda_2$.
\medskip

Let us divide $B$ into a few numbered parts, as described in figure 3.7.
The following short list explains some essential new parts.

\begin{enumerate}
\item The green dashed line, denoted $L$, connects the point $0$ to the midpoint of the arc of $\lambda_1$ bounded by twisted bands $\beta$ and $\gamma$.

\item The three light blue lines bound the $0$-disk.
\end{enumerate}

\begin{center}
\begin{overpic}[scale=0.1]{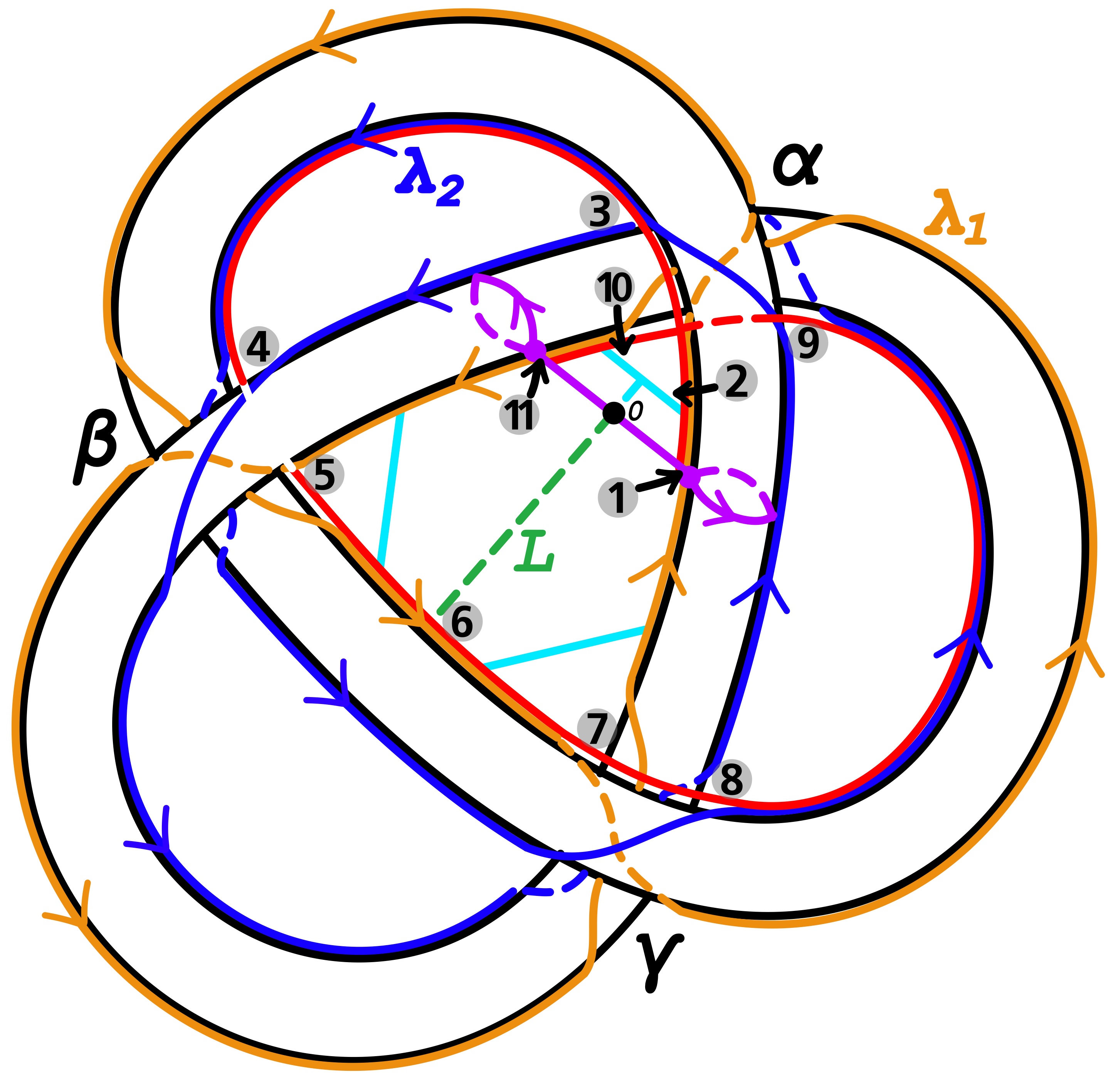}
\end{overpic}
\captionof{figure}{Dividing the boundary into parts}
\end{center}

The following table organizes the different parts of the path along the set $B$. We set our starting point to be $p_2$, denoted by $1$, and our ending point to be $p_1$, denoted by $11$.
\medskip

\begingroup
\renewcommand*{\arraystretch}{1.5} 
\renewcommand*{\baselinestretch}{1.5}
\begin{center}
    \begin{tabular}[c]{ || c | c | c ||}
      \hline  
     Step no. & From point $...$ to point $...$ & Description \\
      \hline \hline 
\rom{1} & $1 \rightarrow 2$ & \makecell{Travel along $\lambda_1$ within the $0$-disk \\ until entering twisted band $\alpha$.} \\
      \hline
\rom{2} & $2 \rightarrow 3$ & \makecell{Enter twisted band $\alpha$ with $\lambda_1$, then travel \\ half a meridian relative to $\lambda_1$  until it \\ coincides with $\lambda_2$. Then exit twisted band $\alpha$.} \\
      \hline
\rom{3} & $3 \rightarrow 4$ & \makecell{Travel along $\lambda_2$ until entering twisted band $\beta$.} \\
      \hline
\rom{4} & $4 \rightarrow 5$ & \makecell{Enter twisted band $\beta$ with $\lambda_2$, then travel \\ half a meridian relative to $\lambda_2$  until it \\ coincides with $\lambda_1$. Then exit twisted band $\beta$.} \\
      \hline
\rom{5} & $5 \rightarrow 7$ & \makecell{Travel along $\lambda_1$ within the $0$-disk while \\ passing through point 6 (intersecting $L$), \\ until entering twisted band $\gamma$.} \\
      \hline
\rom{6} & $7 \rightarrow 8$ & \makecell{Enter twisted band $\gamma$ with $\lambda_1$, then travel \\ half a meridian relative to $\lambda_1$ until it \\ coincides with $\lambda_2$. Then exit twisted band $\gamma$.} \\
      \hline
\rom{7} & $8 \rightarrow 9$ & \makecell{Travel along $\lambda_2$ until entering twisted band $\alpha$.} \\
      \hline
\rom{8} & $9 \rightarrow 10$ & \makecell{Enter twisted band $\alpha$ with $\lambda_2$, then travel \\ half a meridian relative to $\lambda_2$  until it \\ coincides with $\lambda_1$. Then exit twisted band $\alpha$.} \\
      \hline
\rom{9} & $10 \rightarrow 11$ & \makecell{Travel along $\lambda_1$ within the $0$-disk \\ until meeting $p_1$.} \\
      \hline
    \end{tabular}
  \end{center} 
\endgroup
\bigskip

Recall that $\mathcal{H}_1$ describes $\Sigma$'s sweep throughout the trefoil complement, hence the path described in the table has a corresponding path in $\mathcal{H}_1$. The following figure describes this path. The points $1-11$ match the points from figure 3.7.

\begin{center}
\begin{overpic}[scale=0.175]{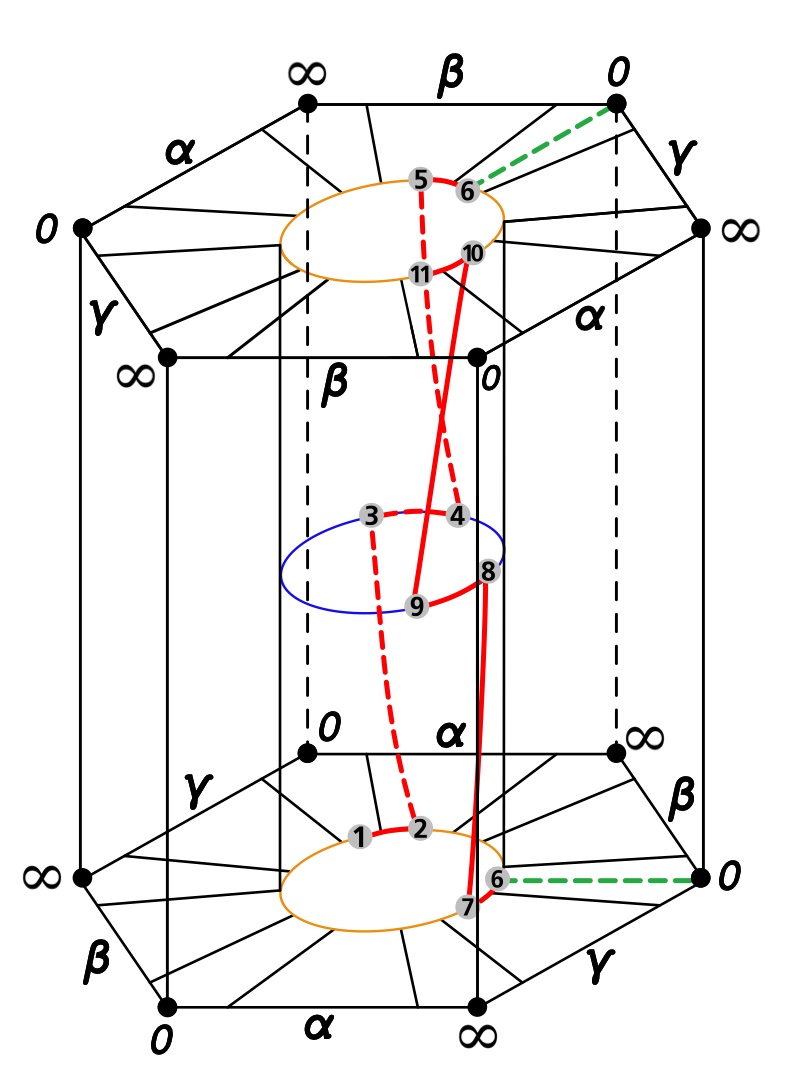}
\end{overpic}
\captionof{figure}{The corresponding path in $\mathcal{H}_1$ }
\end{center}

In conclusion, the entire boundary of the modular template is composed of the path along the set $B$ we have described and of the path along the set $A$, which are connected through the points $p_1$ and $p_2$.

\begin{center}
\begin{overpic}[scale=0.175]{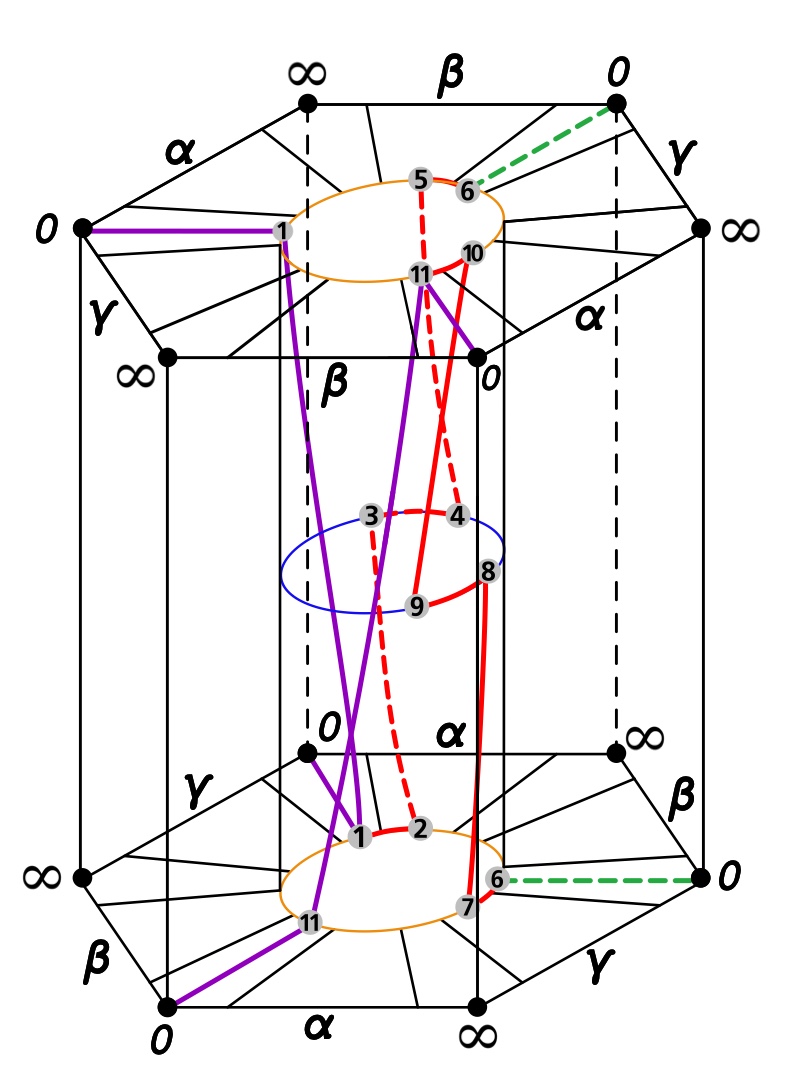}
\end{overpic}
\captionof{figure}{The modular template's boundary in $\mathcal{H}_1$ }
\end{center}

Now, we are ready to embed the template. By gluing identified edges and vertices, the path we described generates a one-component continuous surface. However, when describing the embedding, we will refer to four individual parts, each bounded by some continuous part of the boundary and by a face of $\cH_1$. The first pair of individual parts we consider appears in figure 3.11. The yellow surface is bounded by one-half of the branchline, by $L$, and by the part of the boundary consisting of the points $1,2,3,4,5$ and $6$.
The pink surface is bounded by the other half of the branchline (see also figure 3.10) and by the part of the boundary consisting of the two copies of the point $11$. The two surfaces are glued to one another through the front and back faces of $\cH_1$. In figure 3.12, we see another pair of individual parts, colored green and blue, and similarly, they are glued to one another through the faces that intersect twisted bands $\alpha$ and $\gamma$.
\medskip

The pink and green parts and the yellow and blue parts meet on two different halves of the branchline, and the yellow and green parts are glued together through the two copies of $L$. This results in a one-component surface, which is the embedding of the modular template in $\cH_1$. \hfill $\square$

\begin{center}
\begin{overpic}[scale=0.035]{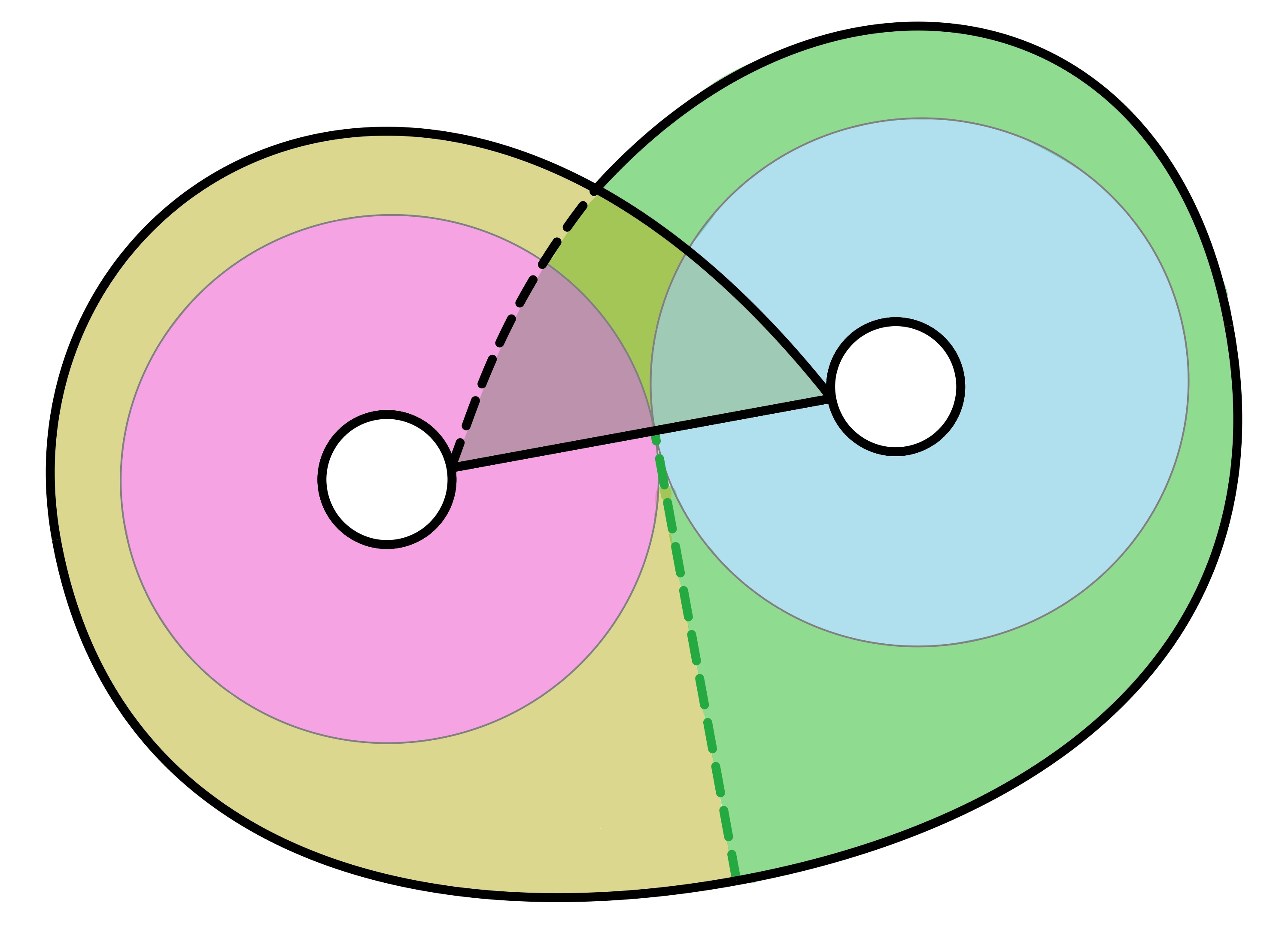}
\end{overpic}
\captionof{figure}{The modular template with color correspondence to figures 3.11 and 3.12}
\end{center}

\begin{center}
\begin{overpic}[scale=0.15]{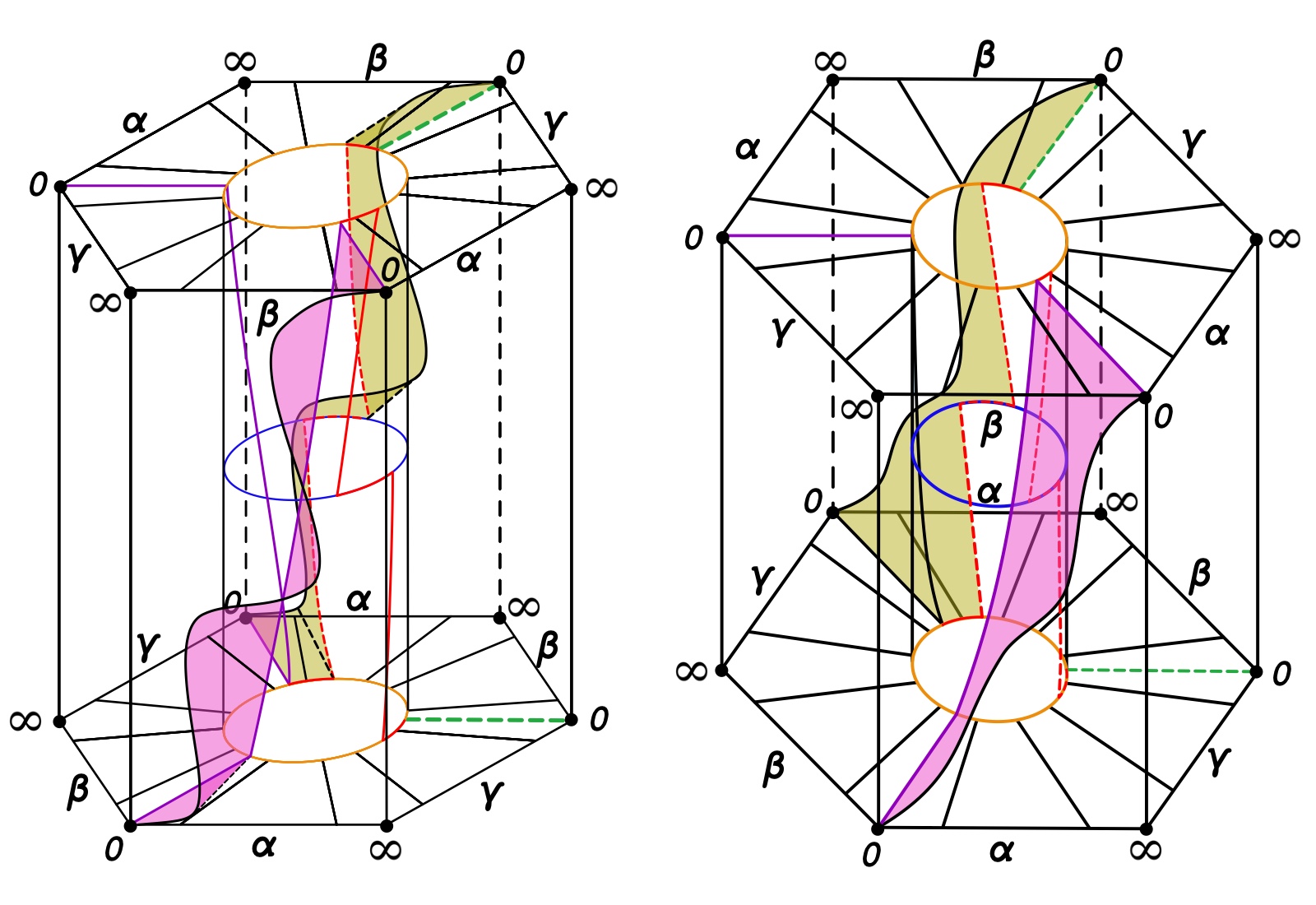}
\end{overpic}
\captionof{figure}{ First pair of individual parts, from two angles}
\end{center}

\begin{center}
\begin{overpic}[scale=0.15]{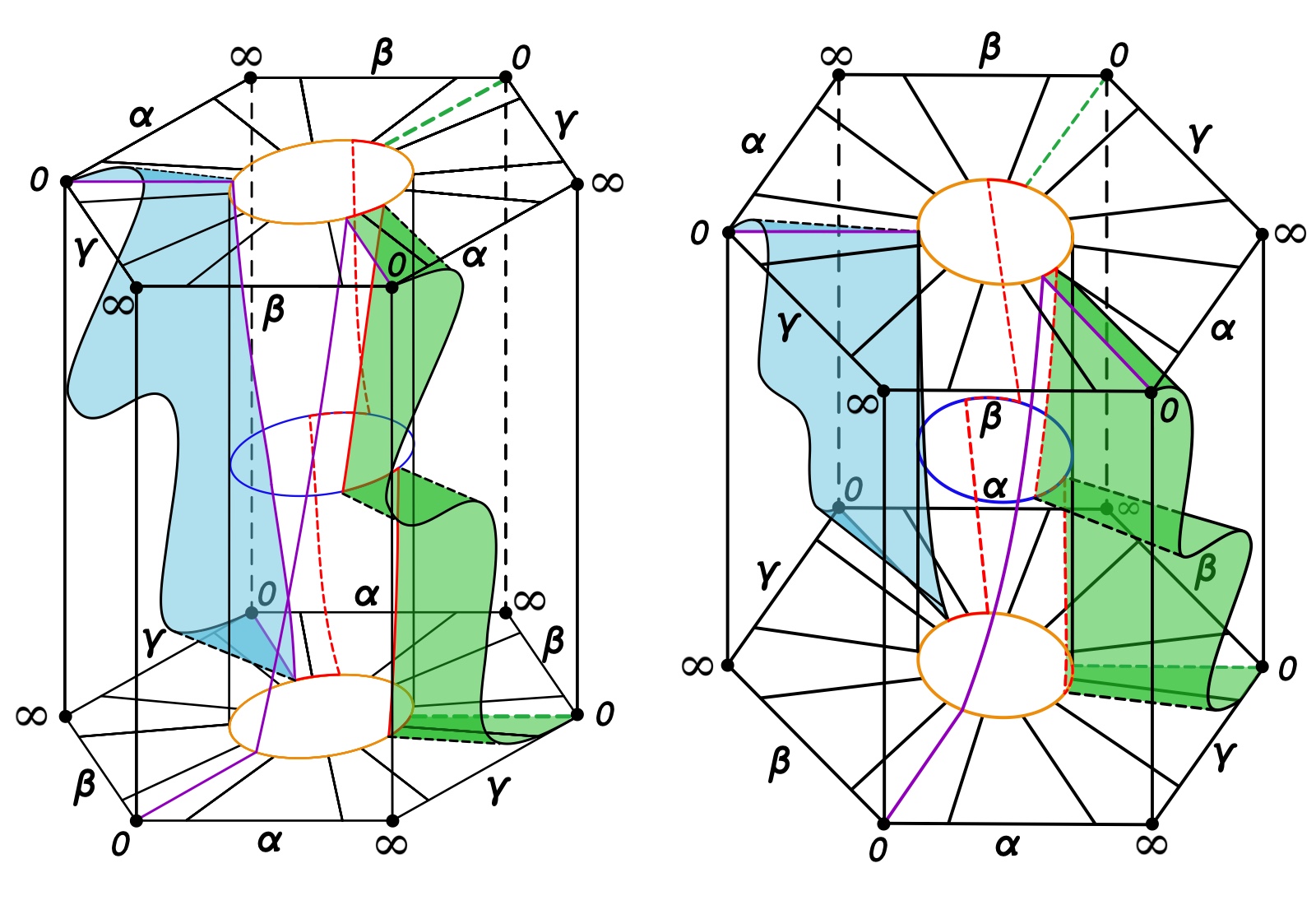}
\end{overpic}
\captionof{figure}{Second pair of individual parts, from two angles}
\end{center}

This gives us the foundation to build a template for every Anosov flow, which is a pre-image of the modular flow.

\theoremm \emph{The template of an Anosov flow, which is the pre-image of the modular flow under $p_n$, is composed of $n$ copies of the modular template, each embedded in a single block of $\cH_n$.}
\medskip

\proof The template of an Anosov flow, which is the pre-image of the map $p_n$, can be constructed in $\cH_n$ by embedding the modular template in each piece $\cH_1$. Then, by stacking $n$ pieces of $\cH_1$ and gluing the $n$th piece to the first with $\varphi$, we obtain the Anosov flow's template. Therefore, by Theorem 3.1, the result follows. \hfill $\square$
\medskip

Theorems 3.1 and 3.2 together prove Theorem 1.1 from the introduction.
\medskip

\remarkk The embedding of the modular template in $\cH_1$, that can be seen in the above figures, is defined up to isotopy. While the surface could be straightened geometrically, we have prioritized a visualization that ensures that the intersections with the faces remain clear, providing a more coherent $3$-dimensional image of the embedding.

\section{Lifting to different covers}

\subsection{Lifting closed geodesics}

The primary purpose of this section is to conclude which periodic orbits lift to closed loops in $\mathcal{H}_n$. Before we state this conclusion, we illustrate a simple example in order to gain intuition.
\medskip

Consider the $R^2L$ geodesic. Figures 4.1.1 (a) and (b) depict the $R^2L$ geodesic in both the modular template and the $\mathcal{H}_1$ model, respectively.

\begin{center}
\begin{overpic}[scale=0.23]{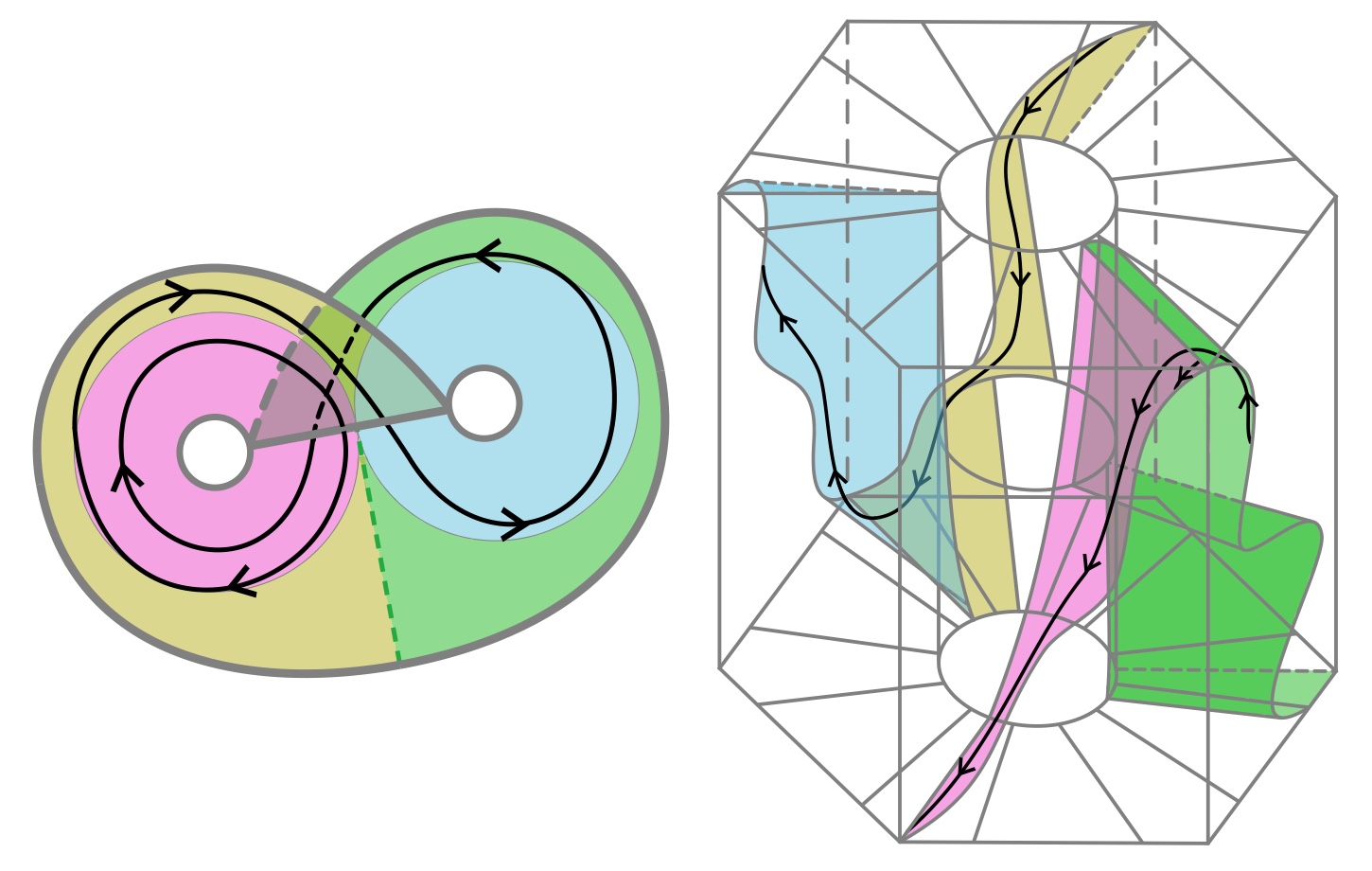}
\put(24.5,-3.5){\color{black}(a)}
\put(73,-3.5){\color{black}(b)}
\end{overpic}
\bigskip

\captionof{figure}{The $R^2L$ geodesic (the colors correspond)}
\end{center}

Along the curve drawn in figure 4.1.1 (b), some parts are directed upwards, and some are directed downwards. The upward-directed movement is consistent with the orientation we chose on the meridian in Note 2.3.5. Thus, every $L$ symbol corresponds to an upward movement along a meridian, while every $R$ symbol corresponds to a downward movement. 
\medskip

However, the lift of the $R^2L$ geodesic to $\cH_n$ results in an open curve. As we begin the movement along the curve, the $R^2$ factor contributes a 2-floor downward movement, while the $L$ factor contributes a 1-floor upward movement. This results in a curve whose endpoints are not identified, hence it is an open curve. This is in contrast to the $RL$ geodesic, which lifts to a curve that travels one floor down and one floor up, thus having the same endpoints and resulting in a closed curve.
\medskip

With that intuition in mind, we now state the main conclusion of this paper, that unfolds how closed geodesics lift to links in $\cH_n$ and how many components they have.
\medskip

\theorem \emph{ Let $K$ be a loop in $S^3\setminus T$, and assume that 
$$K = R^{a_1} L^{b_1} R^{a_2} L^{b_2} \cdot ... \cdot R^{a_t} L^{b_s}, \quad a_i, b_j \in \N.$$ We denote $c:=\gcd(m-\ell, n)$, and write $m = \underset{i}{\sum} a_i$ and $\ell = \underset{j}{\sum} b_j$. Then, one of the following cases occurs:}
\begin{enumerate}
    \item If $c=n$, then $p_n^{-1}(K)$ is a link of $n$ closed components, and each component is encoded exactly as $K$.
    \item If $c=1$, then $p_n^{-1}(K)$ is a one-component link, and $p_n^{-1}(K) = K^n$. Namely, $n$ concatenations of $K$.
    \item If $c>1$, then $p_n^{-1}(K)$ is a link of $c$ closed components, and each component is encoded as $K^{\frac{n}{c}}$.
\end{enumerate}
\medskip

\notation \emph{For the purpose of the following definitions, we introduce a coordinate system on $\cH_n$. Recall $\cH_1 \cong\Sigma\times I$, so we shall denote $\cH_n \cong \faktor{\pare{\overset{n}{\underset{i=1}{\bigcup}} \Sigma\times I_i}}{\sim}$ where $I_i = [i-1,i]$ and the identification is $(x,0) \sim (\varphi(x),n)$. We say $i$ is the floor number.}
\medskip

\deff Let $a = (x_a,t_a)$ and $b = (x_b,t_b)$ be points in $\cH_n$, and denote by $a_0 = (x_a,0)$ and $b_0=(x_b,0)$ their projections to the bottom face of the first floor of $\cH_n$, respectively. Denote the center of the circle $\lambda_1$ by $o$. Then we define the \emph{angle difference between $a$ and $b$} to be $\angle (a,b) :={[\angle a_0ob_0]}_{2\pi}$ and the \emph{floor difference between $a$ and $b$} to be ${[t_a-t_b]}_{n}$.

\lemma \emph{Given an integer $d$, let $P$ be a curve in $\cH_n$ which travels $d$ floors with movement directed either only upwards ($d>0$), or only downwards ($d<0$). Then the angle difference between its endpoints is $r \cdot \frac{2\pi}{6}$, where $d=qn+r$ with $q\in \Z$ and $n > r\in \N\cup \brac{0}$. Moreover, since a curve in $\cH_n$ closes only if the floor and angle difference between its endpoints is zero, $P$ closes to a loop if ${[d]}_{n}=0$.}
\medskip

\proof WLOG we can assume $P$ starts from the bottom of the first floor. From Note 2.3.5, with every floor $P$ travels upwards, the angle difference increases by $\frac{2\pi}{6}$. However, since the $n$th and first floor are glued with $\varphi$, which is a rotation by $\frac{2\pi}{6}$ counter-clockwise, the angle difference decreases by $\frac{2\pi}{6}$ just by identifying the top of the $n$th floor with the bottom of the first. Thus,  traveling all $n$ floors results in a $(n-1) \cdot \frac{2\pi}{6}$ angle difference, but since $n=6k+1$, it results in no angle difference. Now, since $d=qn+r$, $P$ travels all $n$ floors $q$ times and then travels $r$ more floors, resulting in a $r \cdot \frac{2\pi}{6}$ angle difference between its endpoints. \hfill $\square$
\medskip

\proof \textbf{(Theorem 4.1.1)} First, we are only interested in the difference between the number of floors $K$ travels up and down, therefore looking at the sum of $x$ factors and the sum of $y$ factors is sufficient.
\medskip

Let $p_n : \cH_n\longrightarrow S^3 \setminus T$ be as before, and let $\brac{K_i}_{i=1}^n$ be the $n$ pre-images of $K$ assembling the lift $p_n^{-1}(K)$. Denote the starting and ending points of the $i$th component $K_i$ by $p_0^{(i)}$ and $p_1^{(i)}$, respectively, where $p_0^{(i)}$ is at the bottom of the $i$th floor. Notice that for every $i$, $p_0^{(i)}$ differs from $p_1^{(i)}$ by $r={[m-\ell]}_{n}$ floors, and due to Lemma 4.1.4, $\angle (p_0^{(i)},p_1^{(i)}) = r \cdot\frac{2\pi}{6}$. Based on this information we can conclude which $K_i$'s connect and form a closed component. We consider the following values of $c$ as stated in the theorem.
\begin{enumerate}
\item If $c=n$, then $r={[m-\ell]}_{n}=0$. It follows directly from Lemma 4.1.4 that each $K_i$ closes to a loop by itself and forms its own component. Consequently, we obtain a link of $n$ components, where each pre-image is an individual copy of $K$, hence encoded identically.

\item If $c=1$, then the floor difference between the endpoints of each $K_i$ is nonzero, hence $K_i$ cannot form a unique closed component. However, notice that $p_1^{(i)}$ is precisely $p_0^{([i+r]_{n})}$, so $K_i$ connects to $K_{[i+r]_{n}}$ through their endpoints. By repeating this process, all the $K_i$'s connect and form a unique loop. They indeed connect because if we start at $K_1$, we will in turn get to all other components as $\brac{[1+f\cdot r]_{n}}_{f=0}^{n-1} = \brac{1,2,..,n}$. Thus, we pass through all the $K_i$'s exactly once. Notice that the floor difference between the joint curve's endpoints would be $[n\cdot r] \equiv_n 0$, so Lemma 4.1.4 guarantees a one-component link. As for the coding of $p_n^{-1}(K)$, since the endpoint of $K_i$ connects to the starting point of $K_{[i+r]_n}$, the encoding of the resulting closed loop is exactly $K^n$.

\item If $c > 1$, then similarly to the previous case, for each $K_i$ the floor difference between its endpoints is nonzero, so $K_i$ cannot form a unique closed component. Again, some pre-images will connect, but since $c > 1$, only $\frac{n}{c}$ pre-images will form a closed loop at a time. That is true because for a specific $K_i$ the set $\brac{[i+f\cdot r]_{n}}_{f=0}^{n-1}$ consists of only $\frac{n}{c}$ distinct elements. Namely, only $\frac{n}{c}$ out of the $n$ pre-images connect to the $i$th component, which is exactly $\brac{[i+jc]_{n}}_{j=0}^{\frac{n}{c}-1}$ (since $c$ divides $m-l$, there exists $\alpha$ such that $\alpha c = m-l$ and thus $[i+fr]_{n} = [i+f(m-l)]_{n} = [i+f \alpha c]_{n} = [i+jc]_{n}$, where we denote $j=f\alpha$). It follows that we have $c$ such sets, which form the $c$ components. Now, the floor difference between each component's endpoints is $\left[ c\cdot \frac{n}{c}\right] \equiv_n 0 $, and by Lemma 4.1.4, this is a  sufficient condition for each component to form a closed loop. Eventually, we obtain a link of $c$ components. Combining the explanations given in the two previous cases, it is clear that each component is encoded as $K^{\frac{n}{c}}$. \hfill $\square$
\end{enumerate}
\bigskip

\subsection{Commensurability of knots in \texorpdfstring{$\cH_1$}%
     {TEXT}}

Recall that the map $p_n: \cH_n \longrightarrow S^3 \setminus T$ denotes the $n$-sheeted cyclic self-cover of the trefoil complement introduced in $\S$2.3. Each $\cH_n$ carries an Anosov flow that lifts the modular flow, and the periodic orbits of these flows give rise to links of the form $K \cup T \subset S^3$. Since the periodic orbits we study are hyperbolic knots in the trefoil complement \cite{foulon2013contact}, there is a natural connection to the notion of commensurability, which has been extensively studied in the context of hyperbolic knots by Friedl \cite{com:12} and Boileau \cite{boil:11}. In this section, we examine how the links $K \cup T$ are related through their cyclic self-covers and clarify in which sense they are virtually commensurable. We begin with the definition of commensurability.

\deff We say two $3$-manifolds, $N_1$ and $N_2$, are commensurable if they have diffeomorphic finite-sheeted covers. In the case of links, we say two knots are commensurable if their complements satisfy the first part of the definition.
\medskip

\theorem \emph{There are infinitely many periodic orbits $K$ on the modular template, such that the link $K \cup T$ in $S^3$ is commensurable to infinitely many other links with two components where one of them is $T$.}
\medskip

\remark \emph{From Definition 4.2.1, it is clear that the trefoil complement and $\cH_n$ are commensurable, thus our link complements satisfy this definition as well via the cyclic self-covers $\cH_n$. We use the term commensurable, although the term virtually commensurable emphasizes better the fact that the equivalence arises through these finite cyclic self-covers rather than directly between the link complements in $S^3$.}
\medskip

\proof \textbf{(Theorem 4.2.2)}  Let $p_n : \cH_n \longrightarrow S^3\setminus T$ be as before. Denote $K=x^my^\ell$, so from Conclusion 4.1.1, ${p_n}^{-1}(K)$ is a one-component link if and only if $\gcd(m-\ell,n)=1$. It is possible to choose infinitely many $m$ and $\ell$ such that $\gcd(m-\ell,n)=1$ and each choice would form a different knot in $S^3\setminus T$. Therefore, we have infinitely many periodic orbits that lift to a knot. Thus, lifting $K\cup T$ in $S^3$ to the $\cH_n$ will result in a link of two components where one is $T$. Since $\cH_n$ is a cyclic cover and is commensurable to the trefoil complement, then from Definition 4.2.2 these infinitely many links are commensurable to one another. \hfill $\square$
\medskip

\conc \emph{Another consequence of our construction concerns the arithmeticity of the resulting link complements. As established in \cite{pinsky-purcell2023}, certain modular knots, most notably the one corresponding to the $RL$ geodesic, have arithmetic complements. Specifically, the $RL$ complement in the trefoil complement is homeomorphic to the Whitehead link complement, a well-known arithmetic manifold. Since arithmeticity is a commensurability invariant, any finite cover of an arithmetic manifold remains arithmetic. Consequently, the lifts of the $RL$ geodesic to the finite cyclic self-covers $\cH_n$ provide an explicit construction of an infinite family of arithmetic link complements within the self-covers of the trefoil complement. More generally, this principle applies to any modular link in the trefoil complement that is shown to be arithmetic: its pre-image under the covering map $p_n$ will yield a corresponding arithmetic link complement in $\mathcal{H}_n$. This result highlights the non-trivial arithmetic structure emerging in the finite covers of the trefoil complement.}
\medskip

This concludes our results for this paper.
\bigskip

\section{Open questions}

The embedding into $S^3$ of the unit tangent bundle of the modular surface (in which the geodesic flow occurs) is just a single choice. There are infinitely many choices of embedding spaces, each corresponds to a choice of Euler number in the classification of Seifert fibered spaces. Euler number $0$ corresponds to $S^3$ and others correspond to certain lens spaces, but all of them are covered by the trefoil complement. Thus, our embedding of the template helps explore which knot properties are invariant under different embeddings. Another question would be what knot properties of closed orbits of the Anosov flow are preserved in different embeddings.
\medskip

This research has also raised a few thoughts about the connection to the arithmeticity of manifolds. A recent theorem in \cite{bader:21} establishes that a finite-volume hyperbolic $n$-manifold is arithmetic if and only if it contains infinitely many maximal totally geodesic surfaces of real dimension at least $2$, or none at all. These surfaces are associated with hidden symmetries of the manifold, i.e. symmetries between finite index covers of the manifold that do not descend to self-isometries of the base manifold. Many known examples of such manifolds are hyperbolic link complements. Since we have explicitly embedded the $RL$ geodesic in our model, we hope that this framework will assist in studying covering spaces of other knot complements, and in identifying hidden symmetries of these manifolds, thereby contributing to the study of arithmeticity.
\bigskip

\printbibliography

\end{document}